\documentclass[11pt]{amsart}

\usepackage{graphicx} 
\usepackage{times,bbm}
\usepackage{epsfig}
\usepackage{subfigure}
\usepackage{amsthm,amsfonts,amsmath,amssymb}
\usepackage{color}

\advance\textheight 0.1cm
\advance\textwidth 0.1cm
\advance\voffset -0.1cm
\advance\hoffset -0.05cm

\makeatletter\@addtoreset {equation}{section}\makeatother

\newtheorem{thm}{Theorem}[section]
\newtheorem{prop}[thm]{Proposition}
\newtheorem{cor}[thm]{Corollary}
\newtheorem{lem}[thm]{Lemma}

\def\C{\mathcal C}
\def\Cc{{\mathcal C}_c^+}
\def\Chi{\mathbbm{1}}
\def\dist{\mathrm{dist}}
\def\EE{E}
\def\eps{\varepsilon}
\def\HH{{\mathbb H}}
\def\I{\mathcal I}
\def\NN{{\mathbb N}}
\def\Per{{\rm Per}\,}
\def\RR{{\mathbb R}}
\def\SS{{\mathbb S}}
\def\supp{{\rm support\,}}
\def\XX{{\mathbb X}}

\begin{document}

\title{Random polarizations}
\author{Almut Burchard and Marc Fortier}
\address{Department of Mathematics, University 
of Toronto, 40 St. George Street, Toronto, Ontario, Canada M5S 2E4}
\email{almut@math.toronto.edu,reitrof.cram@gmail.com}

\date{November 30, 2012}

\begin{abstract} 
We derive conditions under which random sequences of polarizations
(two-point symmetrizations) 
converge almost surely to the symmetric decreasing rearrangement. 
The parameters for the polarizations are independent random 
variables whose distributions need not be uniform.  
The proof of convergence hinges on an
estimate for the expected distance from the limit that 
yields a bound on the rate of convergence.  In the special case
of i.i.d.~sequences, almost sure convergence holds
even for polarizations chosen at random from suitable small sets.
As corollaries, we find bounds on the rate 
of convergence of Steiner symmetrizations
that require no convexity assumptions, and show
that full rotational symmetry can be achieved by randomly alternating
Steiner symmetrizations in a finite number of
directions that satisfy an explicit non-degeneracy condition.
We also present some negative results on the rate of 
convergence and give examples where convergence fails.
\end{abstract} 
\maketitle

\section{Introduction} 

Many classical 
geometric inequalities were proved by first 
establishing the inequality for a simple geometric transformation,
such as Steiner symmetrization or polarization.
Steiner symmetrization is a volume-preserving rearrangement 
that introduces a reflection symmetry, and polarization 
pushes mass across a hyperplane 
towards the origin. (Proper definitions 
will be given below).  To mention just a few
examples, there are proofs of the isoperimetric 
inequality and Santal\'o's inequality based on the facts that 
Steiner symmetrization decreases perimeter~\cite{Steiner,CS}
and increases the Mahler product~\cite{MP}.
Inequalities for capacities and path 
integrals follow from the observation that 
polarization increases convolution 
functionals~\cite{Wolontis, Dubinin, Ahlfors, BT}
and related multiple integrals~\cite{Schmucki,Mor,Peres-Sousi}.
This approach reduces the geometric inequalities to 
one-dimensional problems (in the case of Steiner symmetrization) or
even to combinatorial identities (in the case 
of polarization). It can also be exploited to characterize equality 
cases~\cite{Ben,Beckner, Schmucki}. A major point
is to construct sequences of the simple rearrangements 
that produce full rotational symmetry in the limit.  

In this paper, we study the convergence of random 
sequences of polarizations to the symmetric decreasing rearrangement.  
The result of $n$ random polarizations of a function $f$
is denoted by $S_{W_1 \dots W_n}f$, where
each $W_i$ is a random variable that determines 
a reflection. We assume that the $W_i$ are independent, but 
not necessarily identically distributed, and derive conditions 
under which 
\begin{equation}
\label{eq:as}
S_{W_1\dots W_n}f \longrightarrow  f^* \quad (n\to\infty) 
\quad \mbox{almost surely}\,.
\end{equation}

Rearrangements have been studied in many different spaces,
with various notions of convergence. We work 
with continuous functions in the topology of uniform convergence, 
while most classical results are stated for compact sets with the
Hausdorff metric. These notions of convergence turn out to be
largely equivalent because of the monotonicity 
properties of rearrangements.

For sequences of Steiner symmetrizations along 
uniformly distributed random directions, convergence is 
well known~\cite{Mani,Volcic}.  It has recently been shown that
certain uniform geometric bounds on the distributions
guarantee convergence for a broad class 
of rearrangements that includes polarization, 
Steiner symmetrization, the Schwarz rounding process, and 
the spherical cap symmetrization~\cite{Schaft1}. 
Among these rearrangements, polarization plays a 
special role, because it is elementary to define, easy to
use, and can approximate the others.
Our conditions for convergence
allow the distribution of the $W_i$ 
to be far from uniform. We also prove
bounds on the rate of convergence, and 
show how convergence can fail. Our results shed new 
light on Steiner symmetrizations.  In particular, we obtain bounds 
on the rate of convergence for Steiner symmetrizations of 
arbitrary compact sets.

\section {Main results} 

Let $\XX$ be either the sphere $\SS^d$,
Euclidean space $\RR^d$, or the standard hyperbolic 
space $\HH^d$, equipped with the uniform Riemannian distance $d(x,y)$,
the Riemannian volume $m(A)$, and a distinguished 
point $o\in\XX$, which we call the origin. 
The ball of radius $\rho$ about a point $x\in\XX$
is denoted by $B_\rho(x)$; if the center is at $x=o$ we 
simply write $B_\rho$. We denote by $\dist(x,A)=\inf_{y\in A} d(x,y)$
the distance between a point and a set, and by
$$
d_H(A,B)=\max\left\{\sup_{x\in A} \dist(x,B), \sup_{x\in B} \dist(x,A)
\right\}
$$
the {\bf Hausdorff distance} between two sets.

If $A$ is a set of finite volume in $\XX$, we denote by $A^*$ the 
open ball centered at the origin 
with $m(A^*)=m(A)$.  We consider nonnegative measurable 
functions $f$ on $\XX$ that vanish weakly at infinity, in the sense that
the level sets $\{x:  f(x)>t\}$ have finite volume for all $t>0$.
(On the sphere, this condition is empty.)
The {\bf symmetric decreasing rearrangement} $f^*$
is the unique lower semicontinuous function that
is radially decreasing about $o$ and equimeasurable  with $f$.  
Its level sets are obtained by replacing the level sets of $f$ with
centered balls,
$$
\{x: f^*(x)>t \} = \{x: f(x)>t \}^* \,.
$$

A {\bf reflection} is an isometry $\sigma$ on $\XX$ with $\sigma^2=I$ 
that exchanges two complementary half-spaces, and has the property 
that $d(x,\sigma y)\ge d(x,y)$ whenever $x$ and $y$ lie in the 
same half-space. On $\SS^d$, we have the reflections at great 
circles, on $\RR^d$ the Euclidean 
reflections at hyperplanes, and in the Poincar\'e ball model of
$\HH^d$ the inversions at $(d\!-\!1)$-dimensional
spheres that intersect the boundary sphere at right angles.
For every point $x\in\XX$ there exists a
$(d\!-\!1)$-dimensional  family of reflections that fix
$x$, and for every pair of distinct points $x,y$ 
there exists a unique reflection that maps $x$ to $y$. 

Let $\sigma$ be a reflection on $\XX$ that does not fix the origin.
For $x\in\XX$, denote by 
$\bar x=\sigma x$ the mirror image of $x$, and let
$$
H^+=\{x: d(x,o)\le d(\bar x, o)\}\,,\quad 
H^-=\{x : d(x,o)\ge d(\bar x,o)\}
$$ 
be the half-spaces exchanged under $\sigma$.
By construction, $o\in H^+$.
The {\bf polarization} of a function $f$ with respect to $\sigma$
is defined by
$$
S f(x)=
\left\{\begin{array}{ll}
   \max\{f(x), f(\bar x)\},  & \text{if } x\in H^+ \,, \\ 
   \min\{f(x), f(\bar x) \}, & \text{if } x\in H^- \,.
\end{array}\right.
$$
For obvious reasons, polarization is also called
{\bf two-point symmetrization}. 

We use a fixed normal coordinate system $x=(r,u)$ 
centered at the origin, where $r=d(x,o)$,
and denote the parameter space by
$\Omega=[0,\infty)\times \SS^{d-1}$.
On $\RR^d$, these are just the standard polar coordinates.
On $\XX=\HH^d$ and $\RR^d$, normal coordinates
define a diffeomorphism from 
$(0,\infty)\times \SS^{d-1}$ to $\XX\setminus\{o\}$,
but on $\XX=\SS^d$ the normal coordinate system degenerates 
at $r=\pi$, where it reaches the south pole.
For $r>0$, let $\sigma_{(r,u)}$ be the reflection
that maps $o$ to  the point with normal coordinates 
$(r,u)$.  The reflections $\{\sigma_{(r,\pm u)}: r > 0\}$ 
generate a one-dimensional group of isometries
of $\XX$.  As $r\to 0$, they converge uniformly to a 
reflection $\sigma_u:=\sigma_{(0,u)}$ that 
fixes the origin and exchanges the half-space 
$H^+_u$ (that has $u$ as its exterior normal at $o$)
with the complementary half-space $H^-_u$.
We do not identify $(0,u)$ with $(0,-u)$ in $\Omega$, although they 
label the same reflection on $\XX$.  If $\omega=(r,u)\in\Omega$ with $r>0$,
the polarization of $f$ with respect to $\sigma_\omega$
is denoted by $S_\omega$.  Given a sequence
$\{\omega_n\}$ in~$\Omega$,
we denote the corresponding sequence of polarizations
by $S_{\omega_1\dots\omega_n}=S_{\omega_n}\circ \dots \circ S_{\omega_1}$.

\begin{figure}[t]
\centering
\subfigure[Polarization 
]{\epsfig{file=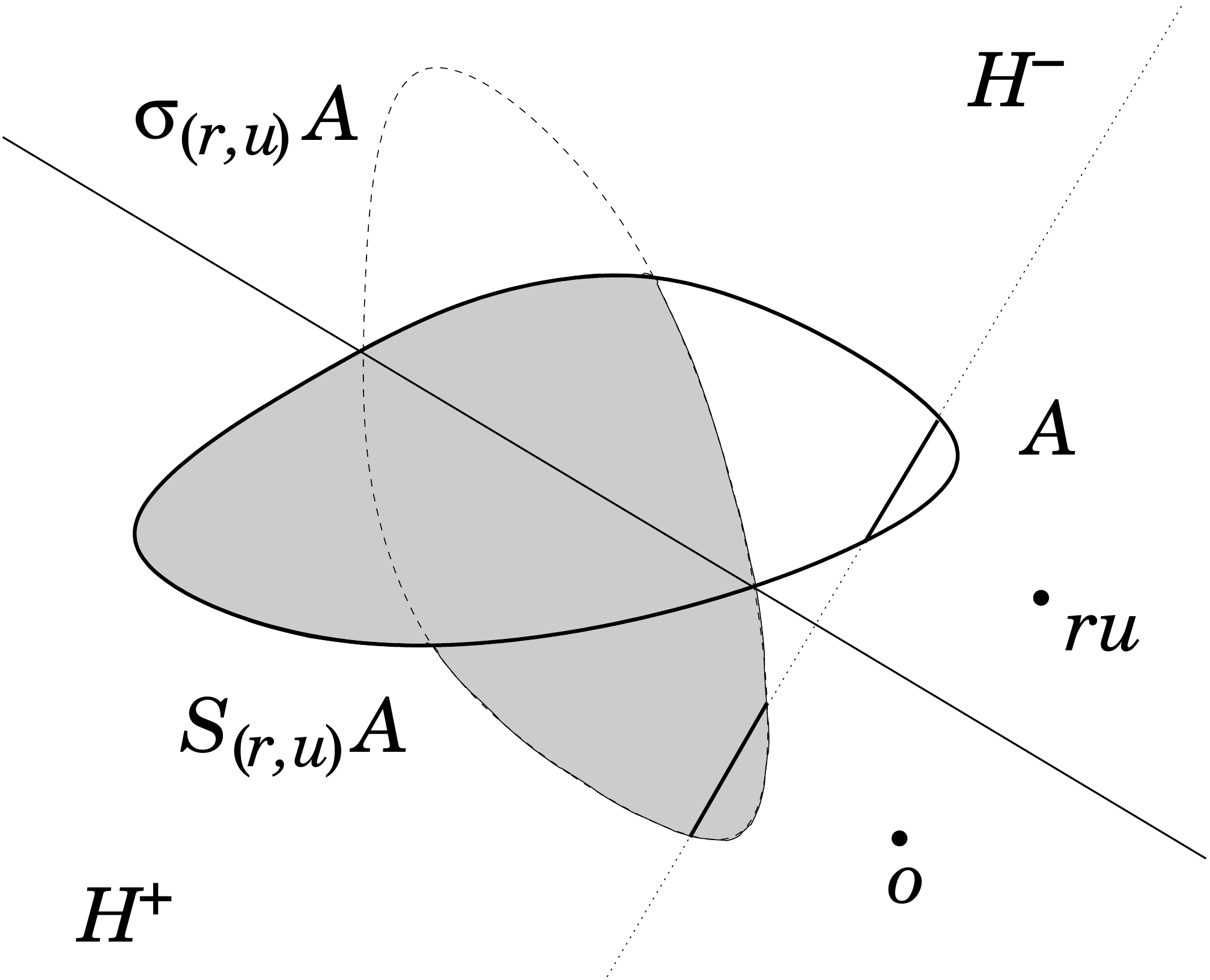, 
width=0.32\linewidth}\label{fig:polarize}
}\hspace{0.12\textwidth}
\subfigure[Steiner symmetrization]{\epsfig{file=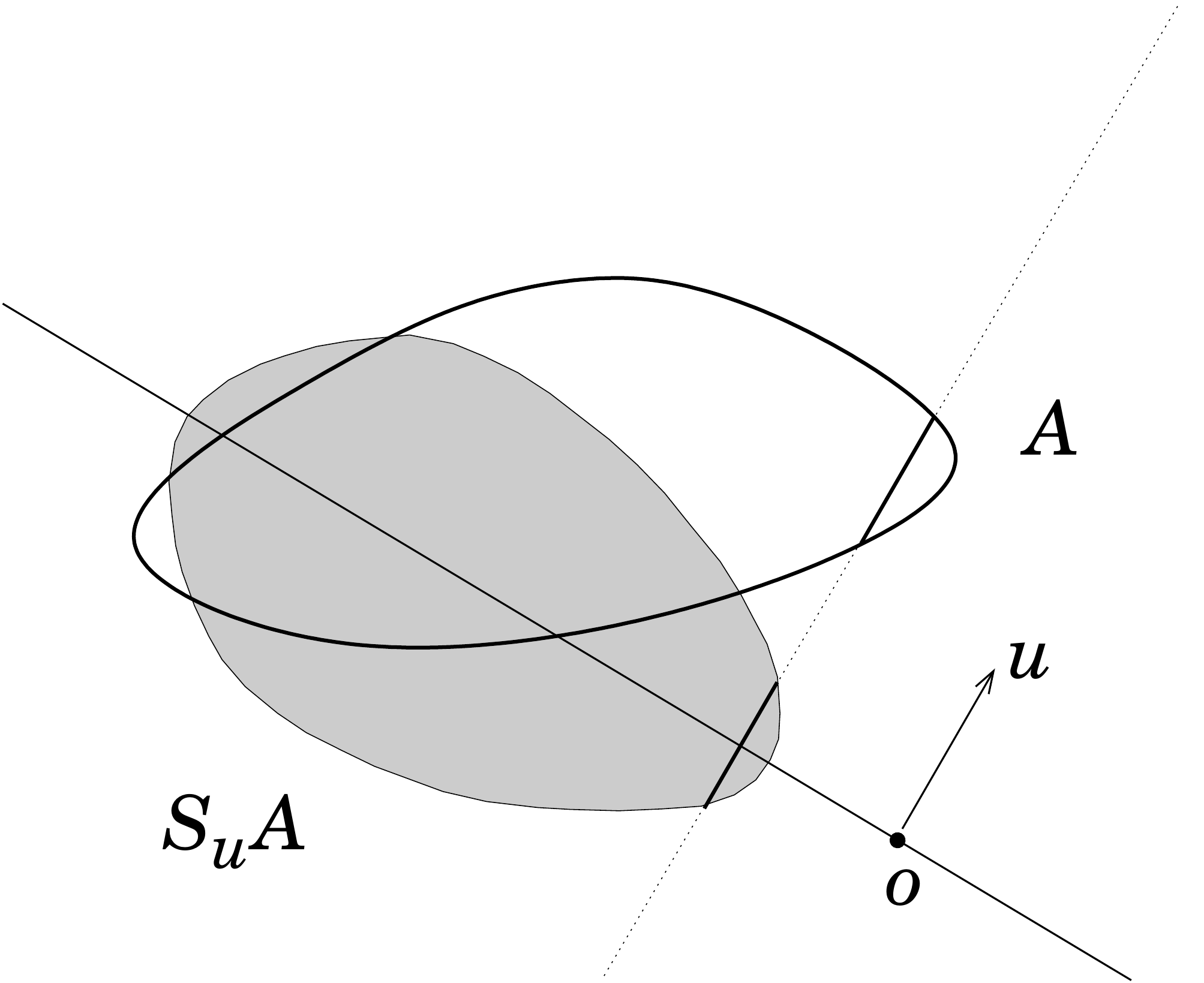, 
width=0.32\linewidth}\label{fig:Steiner}}
\caption{Two simple rearrangements of a set $A$. 
Polarization 
replaces a certain piece
of $A$ in $H^-$ by its reflection in $H^+$. 
Perimeter is preserved but
convexity, smoothness, and non-trivial symmetries can be lost.
Steiner symmetrization
replaces the cross sections of $A$ in a given
direction by centered line segments. This creates
a hyperplane of symmetry and decreases perimeter.
}
\end{figure}

Let $u$ be a unit vector in $\RR^d$, and let $f$
be a nonnegative measurable function that vanishes
at infinity.  The {\bf Steiner symmetrization} 
in the direction of $u$ replaces the restriction of $f$ to each line 
$\{x= \xi+tu : t\in\RR\}$, where $\xi\perp u$,
with its (one-dimensional) symmetric decreasing 
rearrangement. If the restriction of $f$ to such a line
is not measurable or does not decay at infinity, we set the Steiner 
symmetrization of $f$ equal to zero on this line.
We denote the Steiner symmetrization of $f$
by $S_{(0,u)}f$, or simply by $S_uf$.  
By construction, $S_u f$ is symmetric under~$\sigma_u$.
Note that Steiner symmetrization
dominates polarization in the sense that
$$
S_{(r,\pm u)}S_{u} = S_{u} S_{(r,\pm u)} = S_{u}
$$
for every direction $u\in\SS^{d-1}$ and all $r>0$ (see Fig.~1).

Polarization and Steiner symmetrization share 
with the symmetric decreasing rearrangement the 
properties that they are monotone ($f\le g$ implies $Sf\le Sg$), 
equimeasurable ($m(\{Sf>t\})=m(\{f>t\})$ for all $t>0$),
and 
$L_p$-contractive ($||Sf-Sg||_p\le ||f-g||_p$)
for all $p\ge 1$. They also preserve or improve the 
{\bf modulus of continuity}, which we define here as
$$
\eta(\rho)=\sup_{d(x,y)\le \rho} |f(x)-f(y)|\,.
$$
The corresponding rearrangements of a set $A\subset\XX$
are defined by rearranging its indicator function
$\Chi_A$.
Conversely, the rearranged function can be recovered 
from its level sets with the {\bf layer-cake principle}, 
$$
f(x)= \int_0^\infty \Chi_{\{f>t\}}(x)\, dt\,,\quad
Sf(x)= \int_0^\infty \Chi_{S\{f>t\}}(x)\, dt\,.
$$ 

Different from standard conventions, we do not automatically 
identify functions that agree almost everywhere. 
We have chosen the symmetric decreasing 
rearrangement of a function to be lower semicontinuous.
In particular, if $A$ is a set of finite volume, then $A^*$ is an 
{\em open} ball. Polarization and Steiner symmetrization both
transform open sets into open sets. Polarization
also transforms closed sets into closed sets, 
but Steiner  symmetrization does not. 
The literature contains a variant of the symmetric
decreasing rearrangement that preserves compactness, 
where $A^*$ is a {\em closed} centered ball if $A$ has positive volume,
$A^*=\{o\}$ if $A$ is a non-empty set of zero volume, 
and $A^*=\emptyset$ if $A=\emptyset$. Steiner symmetrization is 
again defined by symmetrizing along a family of parallel lines.

A {\bf random polarization} $S_W$ is given 
by a Borel probability measure $\mu$ on 
$\Omega=[0,\infty)\times \SS^{d-1}$
that determines the distribution of the 
random variable $W=(R,U)$, viewed
as the identity map on~$\Omega$.  We assume that
$\mu(R=0)=0$; for $\XX=\SS^d$ we also assume
that $\mu(R>\pi)=0$.  A {\bf random Steiner symmetrization}
$S_U$ is given  by a Borel probability measure on $\SS^{d-1}$,
or equivalently, by a measure on $\Omega$ with
$\mu(R=0)=1$. For sequences of random 
rearrangements $\{S_{W_1\dots W_n}\}_{n\ge 1}$ with each $W_i$ 
independent and distributed according to a measure $\mu_i$ on $\Omega$,
we use as the probability space the infinite product 
$\Omega^{\NN}$ with the product topology,
and with the product measure defined by
$$ P(W_1\in A_1,\ldots, W_n\in A_n)
=\prod_{i=1}^n \mu_i(A_i)\,.
$$
In this view, $W_i=(R_i,U_i)$ is the $i$-th 
coordinate projection on $\Omega^\NN$.

Let $\Cc(\XX)$ be the space of nonnegative
continuous functions with compact support
in $\XX$.  (If $\XX=\SS^d$, this agrees with
the space of all nonnegative continuous functions on
$\SS^d$).  Our first theorem provides a sufficient 
condition for the almost sure convergence of a random sequence 
of polarizations to the symmetric decreasing rearrangement.

\begin{thm}[Convergence of random polarizations]
\label{thm:46}
Let $\{S_{W_1\dots W_n}\}_{n\ge 1}$ be a 
sequence of polarizations on $\XX=\SS^d$, $\RR^d$, or $\HH^d$, 
defined by a sequence of independent random variables $\{W_i\}_{i\ge 1}$ 
on $\Omega$.  If
\begin{equation}\label{assumption-46a}
\sum_{i=1}^{\infty}P(d(\sigma_{W_i}a_i,b_i)<\rho)=\infty
\end{equation} 
for every radius $\rho>0$ and every pair of bounded sequences 
$\{a_i\}$, $\{b_i\}$ in $\XX$ 
with $d(b_i,o)\ge d(a_i,o)+ 2\rho$, 
then
\begin{equation}\label{conclusion-46a}
P\left(\,\lim_{n\rightarrow \infty}\|S_{W_1\dots W_n}f-f^*\|_{\infty}=0 
   \quad \forall f\in \Cc(\XX) \right)=1\,.
\end{equation} 
\end{thm}

At first sight, the conclusion in Eq.~(\ref{conclusion-46a}),
that the random sequence almost surely drives
all functions in $\Cc(\XX)$ {\em simultaneously}
to their symmetric decreasing rearrangements,
looks stronger than Eq.~\eqref{eq:as}.
As we show in the proof of Theorem~\ref{thm:46},
the statements are equivalent, because $\Cc(\XX)$
is separable and polarization contracts uniform distances.
Let $L_p^+(\XX)$ be the space
of nonnegative $p$-integrable functions.
Since polarization also contracts
$L_p$-distances and $\Cc(\XX)$ is dense in $L_p^+(\XX)$,
Eq.~(\ref{conclusion-46a}) extends to
\begin{equation}
\label{eq:Lp}
P\left(\,\lim_{n\rightarrow \infty}\|S_{W_1\dots W_n}f-f^*\|_{p}=0 
   \quad \forall f\in L_p^+(\XX) \right)=1 \qquad (1\le p<\infty)\,.
\end{equation}

The assumption in Eq.~(\ref{assumption-46a}) implies
that infinitely many of the $\mu_i$ assign strictly positive 
measure to every non-empty open set in $\Omega$.
The measures may concentrate 
or converge weakly to zero as $i\to\infty$, but 
not too rapidly. This causes typical random sequences 
to be dense in~$\Omega$. 
We are convinced that
almost sure convergence holds under much 
weaker assumptions on the distribution of the random variables
than Eq.~(\ref{assumption-46a}).
A related question concerns the conditions 
for convergence of non-random sequences $\{\omega_i\}$ 
in $\Omega$.  Clearly, convergence can fail if a sequence 
of polarizations concentrates on a subset of $\Omega$ that 
is too small to generate full rotational symmetry. 
Since the polarization $S_{(r,u)}$ leaves subsets 
of $B_{r/2}$ unchanged, a sequence of
reflections must accumulate near $r=0$ 
to ensure convergence.
 
It is, however, neither sufficient nor necessary 
that the sequence be dense in $\Omega$: on the 
one hand, any given sequence of polarizations can appear as a 
subsequence of one for which convergence fails 
(Proposition~\ref{prop:lower-P}b);
on the other hand, a sequence of polarizations chosen
at random from certain small sets can
converge to the symmetric decreasing
rearrangement (Theorem~\ref{thm:iid}).
Rather, convergence depends on the
ergodic properties of the corresponding reflections
in the orthogonal group $O(d)$.

To state the result, we introduce some more notation.
For $u\in\SS^{d-1}$, let $\tau_u$ be the map 
from $\XX$ to itself that fixes 
the half-space $H_u^+$ 
and reflects the complementary half-space $H_u^-$ by $\sigma_u$.
We visualize $\tau_u$ as folding each
centered sphere down into the hemisphere antipodal to $u$ 
(see Fig.~2a). Given $x\in \XX$ and $G\subset\SS^{d-1}$, we 
refer to the set
$$
{\mathcal O}_{G,x}= \{ \tau_{u_n}\dots \tau_{u_1} x: 
n\ge 0, u_1,\dots, u_n\in G\}
$$
as the {\bf orbit} of $x$ under $G$.

\begin{figure}[t]
\centering
\subfigure[Folding a sphere]
{\epsfig{file=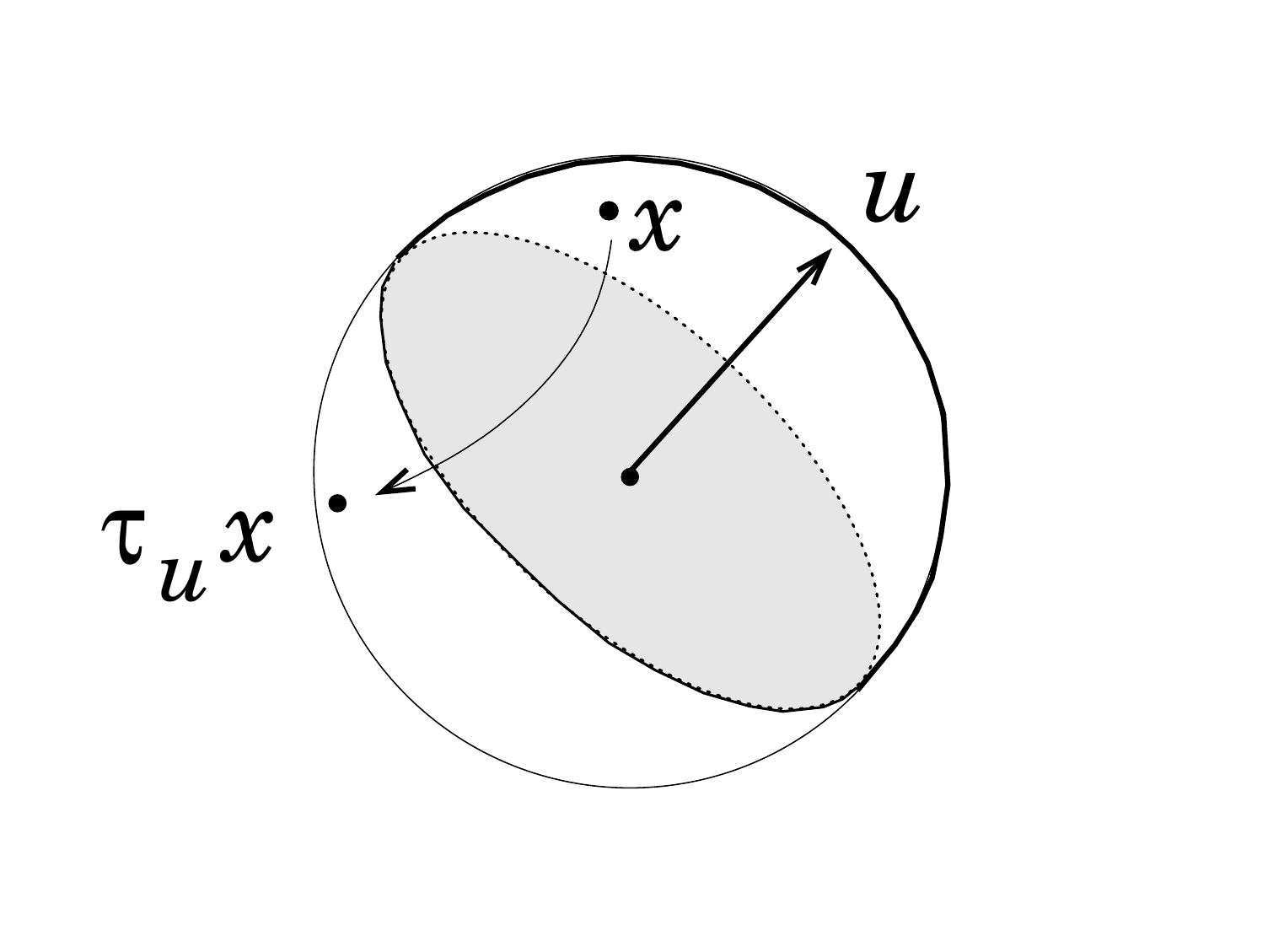, 
width=0.28\linewidth}\label{fig:folding-d}}\hspace{0.12\textwidth}
\subfigure[$G=\{u,v,w\}\subset\SS^1$]{\epsfig{file=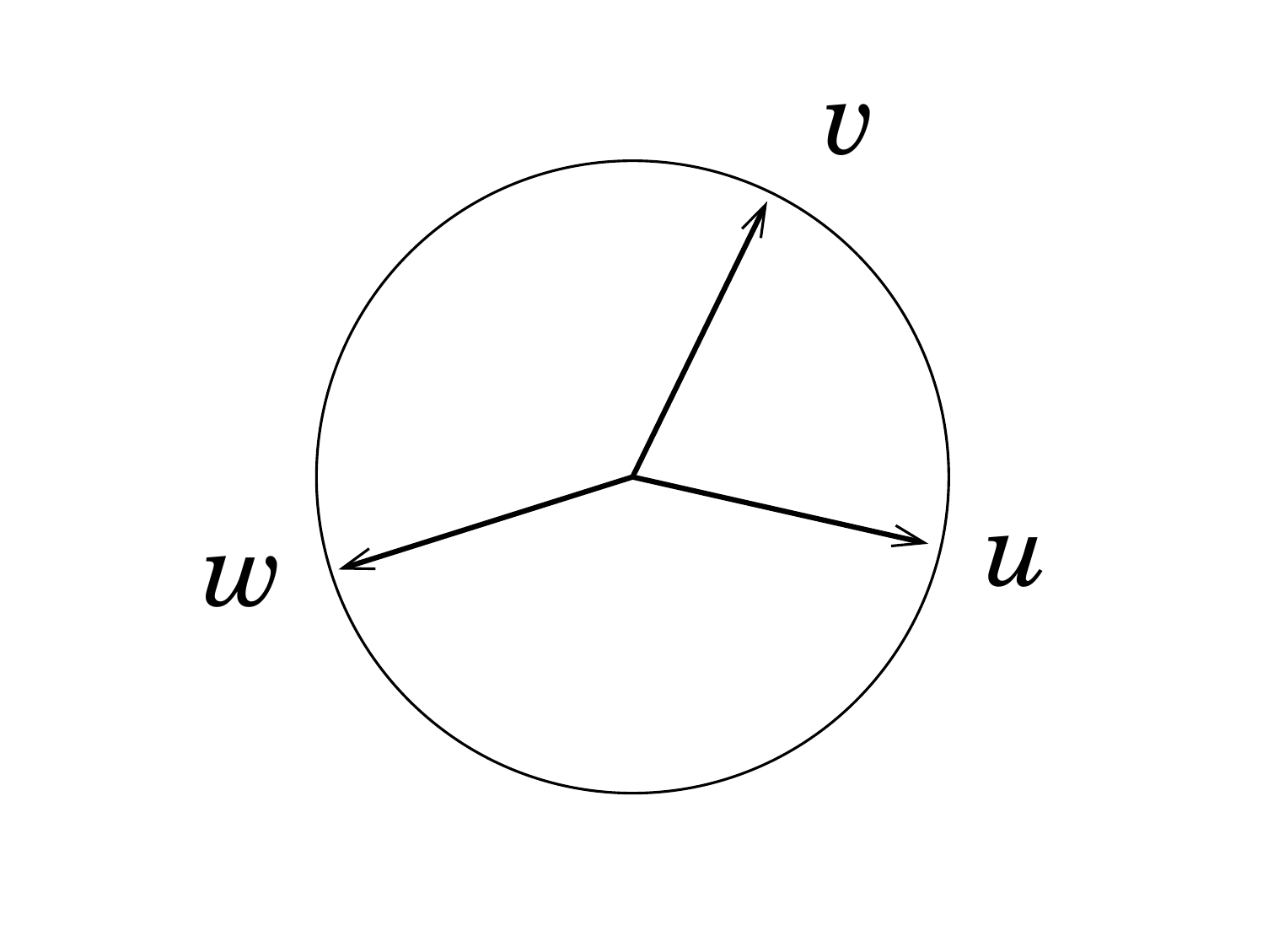, 
width=0.28\linewidth}\label{fig:folding-2}}
\caption{The map $\tau_u$ folds each centered sphere 
in $\RR^d$ into the hemisphere opposite to $u$
across the hyperplane $u^\perp$.
In $d=2$ dimensions, if $u,v,w$ are not contained
in a semicircle and 
enclose angles that are incommensurable with $\pi$,
does $G$ have dense orbits $\mathcal{O}_{G,x}$
in $\SS^1$? Do
$\tau_u,\tau_v,\tau_w$ generate full rotational symmetry?
}
\end{figure}

\begin{thm}[Convergence of i.i.d.~polarizations]
\label{thm:iid}
Let $\{S_{W_1\dots W_n}\}_{n\ge 1}$ be a random sequence of polarizations 
on $\XX=\SS^d$, $\RR^d$, or $\HH^d$, defined
by independent random variables $W_i$  
that are identically distributed according to a 
probability measure $\mu$ on $\Omega$ with $\mu(\{R=0\})=0$.
Let $\supp \mu$ be the smallest closed set 
of full $\mu$-measure in $\Omega$, and set
\begin{equation} \label{eq:iid}
G=\left \{u\in \SS^{d-1}: (0,u)\in \supp \mu\right\}\,.
\end{equation}
If the orbit ${\mathcal O}_{G,x}$ is dense in $\SS^{d-1}$ 
for each $x\in\SS^{d-1}$, then $S_{W_1\dots W_n}$ converges 
to the symmetric decreasing rearrangement
and Eq.~{\rm (\ref{conclusion-46a})} holds.
\end{thm}

In one dimension, polarizations need to accumulate on
both sides of the origin to produce the desired reflection symmetry.
In dimension $d>1$, the precise characterization of subsets 
$G\subset\SS^{d-1}$ that have dense orbits in $\SS^{d-1}$
is an open problem. A necessary condition is
that $G$ be a {\bf generating set of directions}
for the orthogonal group, in the sense that
the finite products $\{\sigma_{u_n}\dots\sigma_{u_1}:
n\ge 0, u_1,\dots, u_n\in G\}$ are dense in $O(d)$.
Also, $G$ cannot be contained in a hemisphere.
A sufficient condition is that
the antipodal pairs $\{u\in G: -u\in G \}$
form a generating set of directions for $O(d)$,
because for every $u\in\SS^{d-1}$
and every $x\in \XX$, either $\sigma_ux=\tau_ux$, or
$\sigma_ux=\tau_{-u}x$.
Must $G$ contain antipodal pairs?
Do $d+1$ directions suffice?  (See Fig.~2b.) 

Generating sets of directions for $O(d)$ are well
understood. For instance, if
{\em (i)} the vectors in $G$ span $\RR^d$;
{\em (ii)} $G$ cannot be partitioned into
two non-empty mutually orthogonal subsets; 
and {\em (iii)} at least one pair of
vectors in $G$ encloses an angle that is not a rational multiple
of $\pi$, then $G$ is a generating set of directions.
(The third condition can be relaxed in dimensions $d\ge 3$.) 
Since $d$ directions $\{u_1,\dots u_d\}$ in general position 
are a generating set for $O(d)$, the hypothesis of Theorem~\ref{thm:iid}
can be satisfied even by measures whose support has
only a finitely many accumulation points.

Theorems~\ref{thm:46} and~\ref{thm:iid}
imply the following statements about
Steiner symmetrization.

\begin{cor} \label{cor:46}
Let $\{S_{U_1\dots U_n}\}_{n\ge 1}$ be a sequence of Steiner 
symmetrizations on $\RR^d$ along independently distributed
random directions $\{U_i\}$ in $\SS^{d-1}$.  

\begin{enumerate} \item [(a)] 
{\rm (Convergence of random Steiner symmetrizations).} 
If \begin{equation}\label{assumption-46b}
\sum_{i=1}^{\infty}P\bigl(d(U_i,v_i)<\rho\bigr)=\infty
\end{equation} 
for every radius $\rho>0$ and every sequence $\{v_i\}$ in $\SS^{d-1}$, then 
\begin{equation}\label{conclusion-46b}
P\left(\,\lim_{n\rightarrow \infty}\|S_{U_1\dots U_n}f-f^*\|_{\infty}=0 
   \quad \forall f\in \Cc(\RR^d) \right)=1\,.
\end{equation}

\item [(b)]
{\rm (Convergence of i.i.d.~Steiner symmetrizations).} 
The same conclusion holds, if, instead,
the random directions $\{U_i\}$ are identically
distributed according to a probability measure $\mu$ on $\SS^{d-1}$
whose support contains a generating set of directions for $O(d)$.
\end{enumerate}
\end{cor}

\section{Related work and outline of the proofs}

The literature contains several different constructions 
for convergent sequences of rearrangements.
In their proof of the isoperimetric inequality,
Carath\'eodory and Study recursively choose
the direction $u_n$ of the next Steiner symmetrization 
such that $A_n=S_{u_n}A_{n-1}$ is as close to 
the ball as possible~\cite{CS}. Lyusternik proposed a sequence 
that alternates Steiner symmetrization
in the $d$-th coordinate direction with Schwarz 
symmetrization in the complementary coordinate hyperplane
and a well-chosen rotation~\cite{Lusternik}.
Brascamp, Lieb, and Luttinger alternate Steiner symmetrization
in all coordinate directions with a rotation~\cite{BLL}.
The constructions of Lyusternik and Brascamp-Lieb-Luttinger 
yield {\em universal} sequences, which work for all nonnegative 
functions on $\RR^d$ that vanish at infinity.

A number of authors have addressed the question of what 
distinguishes convergent sequences of Steiner symmetrizations, 
and how to describe their limits.  Eggleston proved that full 
rotational symmetry can be achieved by iterating Steiner 
symmetrization in $d$ directions that satisfy a non-degeneracy 
condition~\cite[p. 98f]{Eggleston}.
Klain recently showed that iterating any finite set of 
Steiner symmetrizations on a convex body results
in a limiting body that is symmetric under the subgroup of $O(d)$ 
generated by the corresponding reflections~\cite{Klain}. 
On the other hand, Steiner symmetrizations along a 
dense set of directions may or may 
not converge to the symmetric decreasing rearrangement, 
depending on the order in which they are executed~\cite{BKLYZ}.
We note in passing that, although the last three results 
are stated for convex sets, the proofs are readily adapted
to functions in $\Cc(\RR^d)$, with the Arzel\`a-Ascoli theorem
providing the requisite compactness in place
of the Blaschke selection theorem.  By choosing 
the measure in Corollary~\ref{cor:46}b
to be supported on a finite generating set of directions, we obtain an
analogue of Eggleston's theorem for random sequences.

Finding even one convergent sequence of polarizations
is more difficult, because it is not enough to iterate
a finite collection of polarizations. Baernstein-Taylor,
Benyamini, and Brock-Solynin argue by compactness that the set of 
functions that can be reached by some finite number 
of polarizations from a function $f$ contains $f^*$ in its 
closure~\cite{BT,Ben,BS}.  The greedy strategy of Carath\'eodory and Study
also works for the case of polarizations.
Both constructions result in sequences that depend 
on the initial function. A universal sequence
was produced by van Schaftingen~\cite{Schaft2}.

In these papers, considerable effort goes into
the construction of convergent (or non-convergent)
sequences that are rather special. 
The question whether a randomly chosen 
sequence converges with probability one
was first raised by Mani-Levitska~\cite{Mani}.
He conjectured that for compact subsets of $\RR^d$, 
a sequence of Steiner symmetrizations in directions chosen 
uniformly at random should converge in Hausdorff distance
to the ball of the same volume, and verified this for convex sets.

The Mani-Levitska conjecture was settled by van 
Schaftingen for a larger class of 
rearrangements that have the same monotonicity, 
volume-preserving, and smoothing properties as 
the symmetric decreasing rearrangement~\cite{Schaft1}. 
We paraphrase his results for the 
case of polarization.
Van Schaftingen proves the convergence statement
in Eq.~(\ref{conclusion-46a}) under the assumption that
the random variables $W_i$ are independent
and their distribution satisfy the uniform bound
\begin{equation}\label{assumption-VS}
\liminf_{n\to\infty} 
P(d(\sigma_{W_n}a,b)<\rho)>0
\end{equation} 
for every $a,b\in\XX$ and every $\rho>0$.  
In the proof, he 
first constructs a 
{\em universal} sequence, that is, a single non-random sequence
$\{\omega_i\}_{i\ge 1}$ in $\Omega$
such that the symmetrizations 
$S_{\omega_1\dots \omega_n}f$ converge uniformly to $f^*$ 
for every $f\in \Cc(\XX)$.
Eq.~(\ref{assumption-VS}) implies
that typical random sequences closely follow
the universal sequence for arbitrarily long finite segments, 
i.e., for every $\rho>0$ and every integer $N\ge 1$,
$$
P(\exists k: 
d(W_{k+n},\omega_n)<\rho  \ \mbox{for} \   n=1,\dots, N) = 1\,.
$$
After taking a countable intersection over $\rho_N=\frac{1}{N}$ and
$N\in \NN$,
Eq.~(\ref{conclusion-46a}) follows with a continuity argument.

The condition in Eq.~(\ref{assumption-VS}) is stronger
than the corresponding assumption
of Theorem~\ref{thm:46}.  To see this,
let $\{a_i\}$, $\{b_i\}$ be a pair of bounded sequences in  
$\XX$, and choose a pair
of subsequences $\{a_{i_k}\}$, $\{b_{i_k}\}$ that converge 
to limits $a$ and $b$.  For $k$ sufficiently large, 
$$
d(\sigma_\omega a,b)<\frac{\rho}{2}
\ \Rightarrow\ 
d(\sigma_\omega a_{i_k},b_{i_k})<\rho \,.
$$
If Eq.~(\ref{assumption-VS}) holds, then
$P(d(\sigma_{W_{i_k}}a_{i_k},b_{i_k})<\rho)$ 
does not converge to zero,
and the series in Eq.~(\ref{assumption-46a}) diverges.  We 
later show examples that satisfy Eq.~(\ref{assumption-46a}) 
but not Eq.~(\ref{assumption-VS}).

Independently, Vol\v{c}i\v{c} has given a direct geometric
proof for the convergence of Steiner symmetrizations
along uniformly distributed random directions~\cite{Volcic}.
His proof is phrased  as a Borel-Cantelli estimate, 
which suggests that pairwise independence of the $W_i$ might 
suffice for convergence (see~\cite[p. 50-51]{Durrett}). 
Upon closer inspection, there is a conditioning
argument where the independence of the $W_i$ comes 
into play.  It is an open question if convergence can be proved
under weaker independence assumptions.

We are not aware of any prior work on rates of convergence 
for polarizations.  There are, however, some very nice results 
regarding rates of convergence for Steiner symmetrizations of 
convex bodies.  Klartag proved that for every convex body 
$K\subset\RR^d$ and every $0<\eps < 1/2$,
there exists a sequence of  
$ n= \left \lceil cd^4\log^2(1/\eps) \right \rceil $
Steiner symmetrizations $u_1,\dots, u_n$ such that
\begin{equation}\label{eq:Klartag}
d_H(\partial S_{u_1\dots u_n} K,\partial K^*)
\le \eps\cdot \mbox{radius}\,(K^*)\,,
\end{equation}
in other words, $(1-\eps)K^*\subset S_{u_1\dots u_n}K\subset (1+\eps)K^*$.
This means that the distance from a ball decays faster than every
polynomial~\cite[Theorem 1.5]{Klartag}.  
Remarkably, $c$ is a numerical constant 
that depends neither on $K$ nor on the 
dimension. The control over 
the dimension builds on the earlier result of 
Klartag and Milman~\cite{KM} that $3d$ 
Steiner symmetrizations suffice to
reduce the ratio between outradius and inradius of
a convex set to a numerical constant.
Around the same time, Bianchi and Gronchi established 
bounds on the rate of convergence in the other 
direction~\cite{BG}. For each $n$ and every dimension $d$, they 
construct centrally symmetric convex bodies in $\RR^d$
whose Hausdorff distance from a ball cannot be decreased
by {\em any} sequence of $n$ successive Steiner symmetrizations. 
Their construction yields a lower bound on the distance 
from a ball for arbitrary infinite sequences of Steiner 
symmetrizations.  Klartag's results have recently been extended to random 
symmetrizations of convex bodies~\cite{CD}.
It is not known whether convergence is in fact 
exponential, and whether Klartag's convergence estimates 
can be generalized to non-convex sets. 

The proofs of Mani-Levitska, van Schaftingen, 
and Vol\v{c}i\v{c} involve a detailed analysis of typical 
sample paths.  Since they rely on compactness and density 
arguments, they do not yield bounds on the rate of convergence. 
In contrast, Bianchi-Gronchi and Klartag use probabilistic methods 
to find non-random sequences with desired properties. 
The construction of Bianchi and Gronchi takes advantage of 
ergodic properties of reflections. Klartag views the rearrangement 
composed of a random rotation followed by 
Steiner symmetrizations in each of the $d$ coordinate directions
as one step of a Markov chain on convex bodies.
He replaces the Steiner symmetrizations by Minkowski symmetrizations 
to obtain a simpler Markov chain, which acts
on the support function of a convex body
as a random orthogonal projection 
in $L_2$. Since this simpler process is a strict contraction 
on the spherical harmonics of each positive order, 
the support function converges
exponentially (in expected $L_2$-distance) 
to a constant. He finally obtains Eq.~(\ref{eq:Klartag})
from a subtle geometric comparison argument. 

We combine an analytical
approach similar to Klartag's with the geometric
techniques used by Vol\v{c}i\v{c}.
The sequence 
$\{S_{W_1\dots W_n}f\}_{n\ge 1}$ defines a Markov chain
on the space $\Cc(\XX)$. 
We use that the functional
\begin{equation} \label{eq:I}
\I(f)=  \int_{\XX} f(x) \, d(x,o)\, dm(x)
\end{equation}
decreases under each polarization, and 
make Vol\v{c}i\v{c}'s conditioning argument explicit
by appealing to the Markov property. 
Here, $dm(x)$ denotes integration with
respect to the standard Riemannian volume on 
$\XX=\SS^d$, $\RR^d$, or $\HH^d$. 
For the proof of Theorem~\ref{thm:46},
we quantify the expected value of the
drop $\I(f)-\I(S_Wf)$ 
in terms of $||f-f^*||_\infty$
and the modulus of continuity of $f$.
Since the expected drop goes to zero,
$S_{W_1\dots W_n}f$ converges uniformly to~$f^*$.

For the case of i.i.d.~polarizations considered in 
Theorem~\ref{thm:iid}, the challenge is that their distribution may 
be supported on a small set. Here, we resort 
to a compactness argument.
By monotonicity, $\I(S_{W_1\dots W_n}f)$ approaches
a limiting value. Under the assumptions of the theorem,
the drop of $\I$ has {\em strictly positive expectation}
unless $f=f^*$ (Lemma~\ref{lem:OG}).
This forces the limits of convergent subsequences to be invariant
under a  family of transformations (the folding maps 
$\tau_u$ parametrized by Eq.~\eqref{eq:iid}),
which play the role of {\em competing symmetries}~\cite{CL}: 
the only functions that are invariant under the entire family 
are constant on each centered sphere.

Our estimates for the expected drop of $\I$
imply bounds on the rate of convergence 
that depend on the modulus of continuity of $f$
and the distribution of the $W_i$. In the case where 
the $W_i$ are uniformly distributed on a suitable subset of
$\Omega$, we show that there exists a 
numerical constant $c$ such that
$$
\EE(||S_{W_1\dots W_n}f-f^*||_\infty)
\le c L\,{\rm Lip}(f)\, n^{-\frac{1}{d+1}}
$$
for every Lipschitz continuous nonnegative function $f$ on
$\RR^d$ with support in $B_L$ (Proposition~\ref{prop:uniform}).
On the other hand, there exist Lipschitz continuous 
functions $f$ with support in $B_L$ such that
$$
\EE(||S_{W_1\dots W_n}f-f^*||_\infty) 
\ge c \, ||f-f^*||_\infty  \, q^n\,,
$$
where $c>0$ and $q\in (0,1)$ are numerical constants
(Proposition~\ref{prop:lower-P}a).

For Steiner symmetrization, we use that
\begin{equation}\label{eq:trump}
\I(f^*)\ \le \ \I(S_uf) 
\ \le \ \I(S_{(r,\pm u)}f)\ \le\ \I(f)
\end{equation}
for every $u\in\SS^{d-1}$ and all $r>0$
to bound the expected value of the drop  $\I(f)-\I(S_Uf)$
under a random Steiner symmetrization from below by
the corresponding estimate for a random polarization 
(Corollary~\ref{cor:46}). By the same token,
the power-law bounds on the rate of convergence
extend to Steiner symmetrizations along uniformly distributed 
directions (Corollary~\ref{cor:uniform-S}). 
Since we ignore that Steiner
symmetrization reduces perimeter, these bounds cannot be sharp, but 
to our knowledge they are the only available bounds 
that do not require convexity.
It is an open question whether the sequence 
converges exponentially, and how the rate of convergence depends 
on the dimension.  Is it more effective to alternate Steiner 
symmetrizations along the coordinate directions with a 
random rotation, as in~\cite{Klartag}?
Does it help to adapt the sequence to the function?
Do polarizations converge more slowly, perhaps following 
a power law?  

\section{ Almost sure convergence}

We start by preparing some tools for the proof of the main results.
Let $\I$ be the functional defined in Eq.~\eqref{eq:I}.
The first lemma is a well-known identity, which is related
to the Hardy-Littlewood inequality $\int fg\le\int f^*g^*$.
We reproduce its proof here for the convenience of the reader.

\begin{lem}[Polarization identity] \label{lem:HL}
Let $f$ be a nonnegative measurable function
with $\I(f)<\infty$, and let $S_\omega$ be  a polarization.
Then 
$$
\I(f)-\I(S_\omega f) = 
\int_{\XX} 
[f(\sigma_\omega x)\!-\!f(x)]^+\,
[d(\sigma_\omega x,o)\!-\!d(x,o)]^+\, dm(x)\,.
$$
In particular, $\I(f)>\I(S_\omega f)$ unless 
$S_\omega f=f$ almost everywhere.
\end{lem}

\begin{proof} We rewrite the functional as an 
integral over the positive half-space $H^+$ associated with $\omega$,
$$
\I(f)-\I(S_\omega f) = \int_{H^+}
\bigl\{
(f(x) \!-\!S_\omega f(x)) d(x,o)
+ (f(\bar x)\!-\!S_\omega f(\bar x))d(\bar x,o)\bigr\}\, dm(x)\,,
$$
where $\bar x=\sigma_\omega(x)$.
If $f(x)\ge f(\bar x)$ for some $x\in H^+$, 
then the values of $S_\omega f$ at $x$ and $\bar x$
agree with the corresponding values of $f$,
and the integrand vanishes at $x$. If, on the other hand
$f(x)<f(\bar x)$, then the values are swapped
for $S_\omega f$, and the integrand becomes
$(f(x)\!-\!f(\bar x))(d(x,o)\!-\!d(\bar x,o))$,
where both factors are negative. We switch the signs,
collect terms, and integrate to obtain the claim.
\end{proof}

The next lemma is the key ingredient in the proof of Theorem~\ref{thm:46}.

\begin{lem}[Expected drop of $\I$] \label{lem:key}
Let $f$ be a nonnegative continuous function with 
compact support in $B_L\subset\XX$ for some $L>0$
and modulus of continuity $\eta$.
Set $\eps=||f-f^*||_\infty$, let 
$\rho>0$ be so small that
$\eta(\rho)\le \frac{\eps}{8}$,
and let $W=(R,U)$ be a random variable on $\Omega$, as described above.
Then 
\begin{equation}
\label{eq:key-pol}
\EE(\I(f)-\I(S_Wf)) \ge C_\eps\cdot
\inf_{x, b} P( d(\sigma_W x, b)<\rho)\,,
\end{equation}
where $C_\eps= \eps\rho\, m(B_\rho)/2$,
and the infimum extends over $x,b$ with
$d(x,o)+2\rho \le d(b,o)\le L-\rho$.
Furthermore, on $\XX=\RR^d$, 
\begin{equation}
\label{eq:key-Steiner}
\EE(\I(f)-\I(S_Uf)) \ge C_\eps'\cdot
\inf_{v\in\SS^{d-1}} P(2L \sin d(U,v)<\rho)\,,
\end{equation}
where $C_\eps'=\eps\rho\, m(B_\rho)/8$.
\end{lem}

\begin{proof} 
We first construct a pair of points $a,b\in\XX$ such that
$$
d(b,o) \ge  d(a,o)+ 4\rho\,,\quad f(b)\ge f(a)+\frac{\eps}{2}
$$
(see Fig.~3a).  By assumption, there exists a point $x_0$ with 
$|f(x_0)-f^*(x_0)|= \eps$. Set $t=\frac12 (f(x_0)+f^*(x_0))$,
let $A=\{ x : f(x) > t\}$,
and let $A^*$ be the  corresponding level set of $f^*$.
If $f(x_0)<f^*(x_0)$, we set $a=x_0$. By construction, 
$a\in A^*\setminus A$. Since this set is open and
non-empty, it has positive volume, and therefore 
$A\setminus A^*$, having the same volume,
is non-empty. Let $b\in A\setminus A^*$. 
Then $f^*(a)-f^*(b) > f^*(a)-t=\eps/2$.
Similarly, if $f(x_0)>f^*(x_0)$, we set $b=x_0\in A\setminus A^*$, 
find $a\in A^*\setminus A$, and note that
$f^*(a)-f^*(b) >t-f^*(b)=\eps/2$.
Since the modulus of continuity of $f$ is valid also
for $f^*$ and $\eta(4\rho) \le \eps/2$, we have
$d(b,o)-d(a,o)\ge 4\rho$.

By Lemma~\ref{lem:HL} and Fubini's theorem, a random polarization
$S_W$ satisfies
\begin{eqnarray*}
\EE(\I(f)-\I(S_Wf)) \hspace{-2cm}&&\\
&=& \EE\left(
\int_{\XX}
[f(\sigma_W x)-f(x)]^+\,
[d(\sigma_W x,o)-d(x,o)]^+\, dm(x)\right)\\
&\ge& \frac{\eps\rho}{2} 
\int_{B_\rho(a)}
P( d(\sigma_W x,b)<\rho)\, dm(x)\,,
\end{eqnarray*}
because the choice of $a$ and~$b$ ensures that
$f(\sigma_W x)-f(x)\ge\eps/4$ and
$d(\sigma_W x,o)-d(x,o)\ge 2\rho$ 
for all $x\in B_\rho(a)$ with $d(\sigma_Wx,b)<\rho$.
Eq.~\eqref{eq:key-pol} follows by minimizing over 
$x$, $a$ and $b$ and evaluating the integral.

\begin{figure}[t]
\centering
\subfigure[Choice of $a$, $b$ in Lemma~\ref{lem:key}]
{\epsfig{file=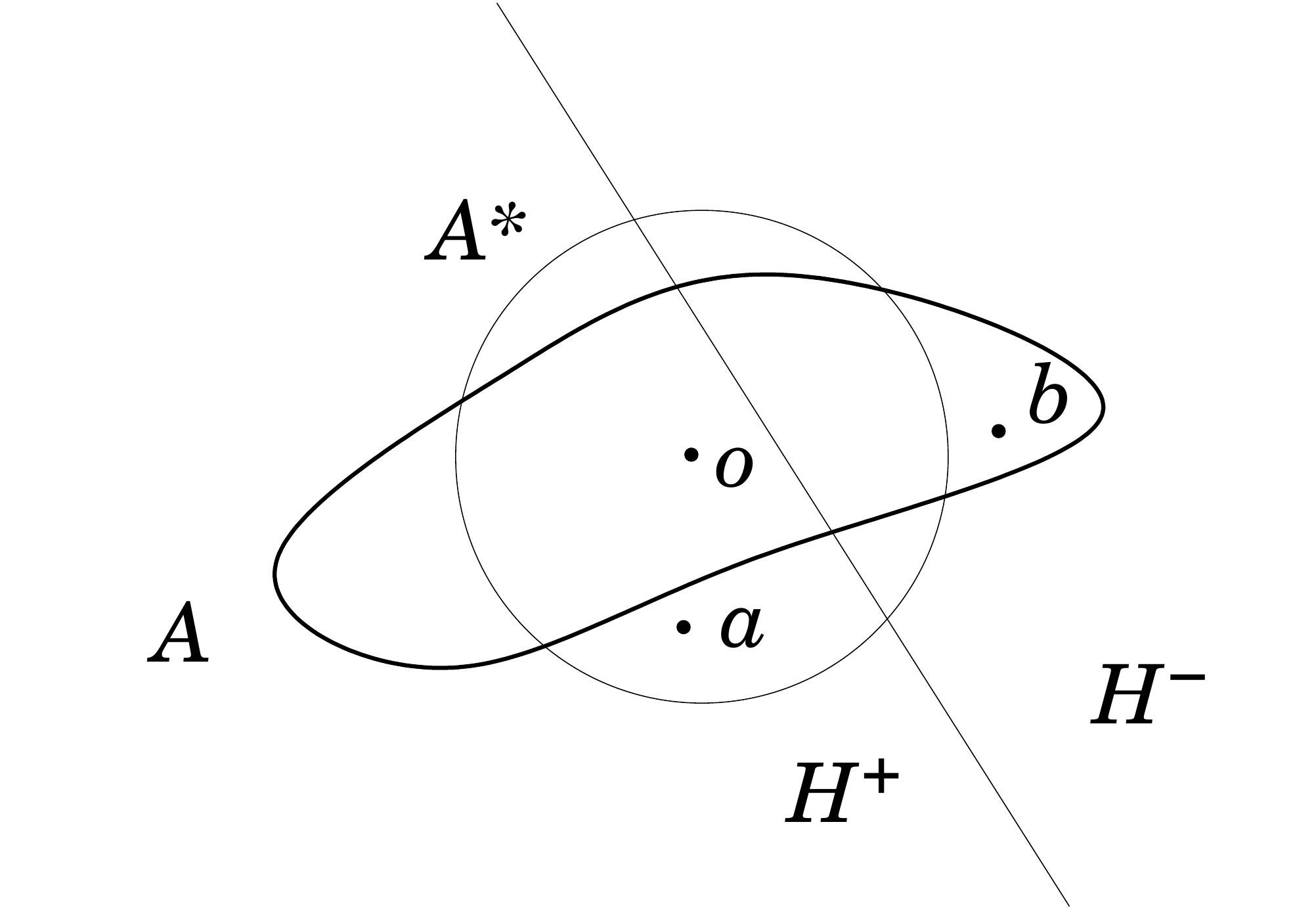, 
width=0.32\linewidth}\label{fig:folding-pair}
}\hspace{0.12\textwidth}
\subfigure[The estimate in Eq.~\eqref{eq:Delta}]{\epsfig{file=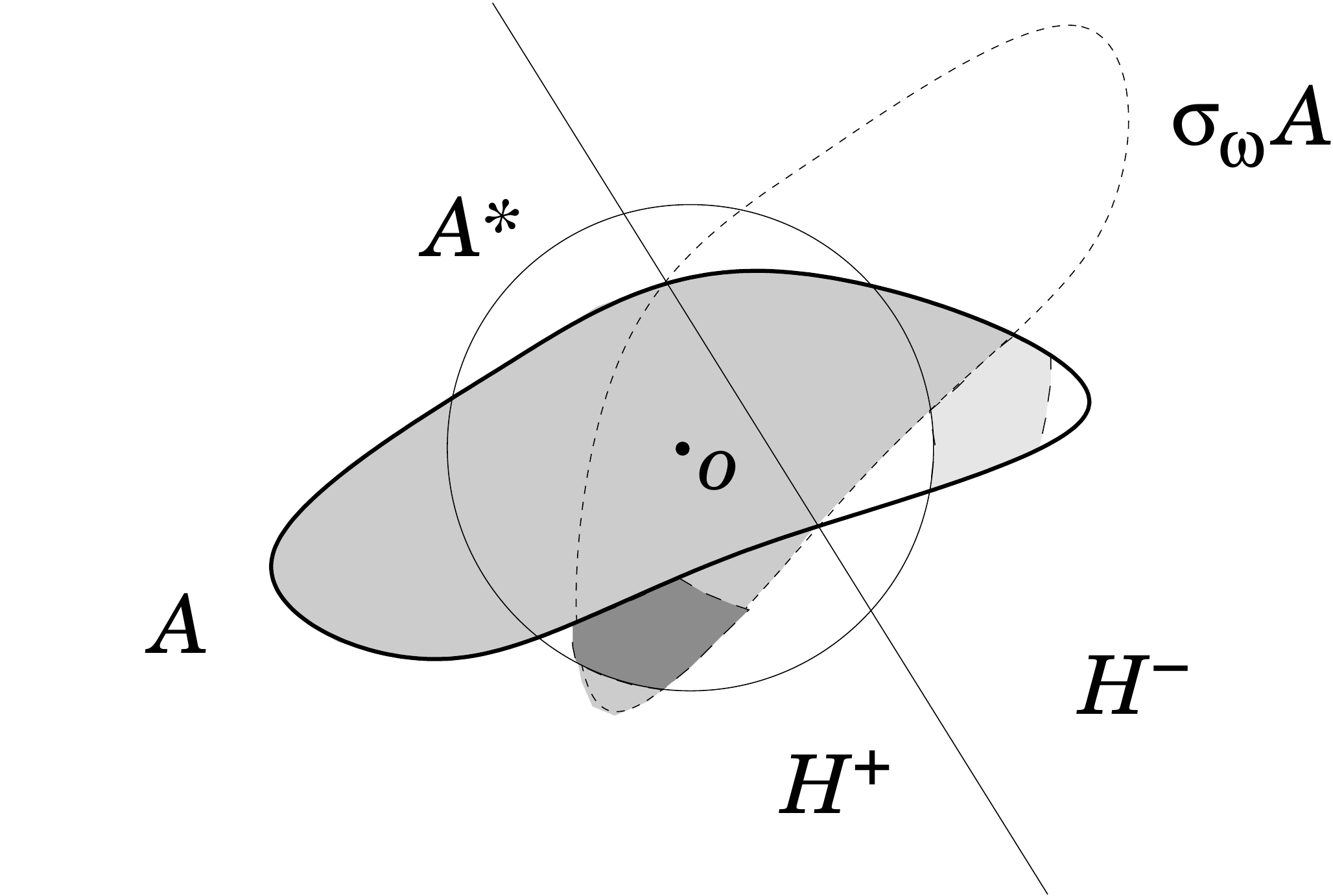, 
width=0.32\linewidth}\label{fig:Delta}}
\caption{Polarization swaps the
part of $A\setminus \sigma_\omega A$ 
that lies in in $H^-$ with its mirror image in $H^+$.
If $A$ is a level set of $f$
and $\sigma_\omega$ is the reflection
that maps $a$ to $b$, then $\I(S_\omega f)<\I(f)$.
The volume of $A\bigtriangleup A^*$ decreases by the 
combined volume of the light and dark shaded subsets.
}
\end{figure}

For a random Steiner symmetrization $S_U$, we use Eq.~(\ref{eq:trump}) 
to obtain
\begin{eqnarray*}
\notag
\EE(\I(f)-\I(S_Uf))
\hspace{-2cm}&&\\
&\ge& \EE\left(
\sup_{r>0,\pm } \ 
\int_{\XX}
[f(\sigma_{(r, \pm U)}x)-f(x)]^+\,
\bigl[\,|\sigma_{(r,\pm U)}x|-|x|\,\bigr]^+\, dm(x)\right)\\
&\ge& \EE\left(
\Chi_{\inf
|\sigma_{(r,\pm U)}a-b|<\rho}
\int_{B_\rho(a)} \frac{\eps\rho}{4}\, dm(x)\right)\\
&\ge& \frac{\eps\rho}{8} \, m(B_\rho(a))
\cdot 
P(2L\sin d(U,v)<\rho)\,,
\end{eqnarray*}
where $v$ is the unit vector in the direction
of $b-a$, and $d(U,v)$ is the enclosed angle.
In the second line, the infimum runs over $(r>0,\pm)$,
and we have used that $f(\sigma x)-f(x)\ge \eps/8$ 
and $|\sigma x| -|a|>\rho$ whenever
$|x-a|<\rho$ and $|\sigma a-b|<\rho$.
In the last line, we have estimated
the infimum by
$$\inf_{r>0,\pm} |\sigma_{(r,\pm U)}a-b|
\ =\ \inf_{t\in\RR} |a+tU-b|
\ \le \ (|a|+|b|)\sin d(U,v)\,,
$$
and applied Fubini's theorem.
\end{proof}

\begin{proof}[Proof of Theorem~\ref{thm:46}]
Given $f\in \Cc(\RR^d)$,  let $F_n=S_{W_1\dots W_n}f$ be the 
result of $n$ random polarizations of $f$.  Since $F_n=S_{W_n}F_{n-1}$,
the sequence $\I(F_n)$ decreases monotonically and satisfies 
$$
\I(f)\ \ge \ \I(F_{n-1}) \ \ge \ \I(F_n) \ \geq  \ \I(f^*)\,.
$$ 
By writing the difference as a telescoping sum and
taking expectations, this implies that 
\begin{eqnarray} 
\notag \I(f)-\I(f^*) \hspace{-2cm}&&\\
\notag
&\ge&  \EE \left( \sum_{n=1}^\infty \I(F_{n-1})-
\I(S_{W_n}F_{n-1})\right)\\
\notag &=& \sum_{n=1}^\infty \ 
\EE( \EE(\I(F_{n-1})-\I(S_{W_n}F_{n-1})
\ \vert\ W_1\dots W_{n-1})) \\
\label{eq:tata}
 &\ge& C_\eps \cdot \sum_{n=1}^\infty 
\left\{
\inf_{x,b} P( d(\sigma_{W_n}x,b)>\rho)
\cdot P( ||F_{n-1}-f^*||_\infty\ge \eps)\right\},
\end{eqnarray}
where the infimum extends over all $x,b$ with
$d(x,o)+2\rho \le d(b,o)\le L-\rho$,
and $C_\eps$, $\rho$, and $L$
are positive constants that depend on $f$.
We have used the Markov property in the second step,
and applied Eq.~(\ref{eq:key-pol}) of Lemma~\ref{lem:key} in the third.
In particular, the sum in Eq.~(\ref{eq:tata})
converges. Since the first factors
in the product are not summable by Eq.~\eqref{assumption-46a},
the second factors must have zero as an accumulation
point. By monotonicity, they converge to zero.
Since $\eps>0$ was arbitrary, we conclude that
$$
P\Bigl(\lim_{n\to\infty} ||F_n-f^*||_\infty= 0 \Bigr)=1\,.
$$
This establishes Eq.~\eqref{eq:as}. 

To complete the proof, we
choose a countable dense subset $\mathcal G \subset \Cc$.
Let $\{\omega_i\}_{i\ge 1}$ be a sequence in $\Omega$.
Since polarizations and the symmetric decreasing rearrangement
contract uniform distances, we have for
every pair of functions
$f,g\in\Cc$ and every $n\ge 1$,
$$
||S_{\omega_1\dots \omega_n}f-f^*||_\infty
\le 2||f-g||_\infty + ||S_{\omega_1\dots \omega_n}g-g^*||_\infty\,.
$$ 
We take $n\to\infty$ and minimize over $g\in\mathcal G$
to obtain, by the density of $\mathcal G$,
\begin{align*}
\lim_{n\to\infty} ||S_{\omega_1\dots \omega_n}f-f^*||_\infty
&\le\inf_{g\in\mathcal G}
\left\{ 2||f-g||_\infty + 
\lim_{n\to\infty}||S_{\omega_1\dots \omega_n}g-g^*||_\infty\right\}\\
&\le \sup_{g\in\mathcal G} \ \lim_{n\to\infty}
||S_{\omega_1\dots \omega_n}g-g^*||_\infty\,.
\end{align*}
Since $\mathcal G$ is countable, it follows that
\begin{eqnarray*}
P\left (\exists f\in\Cc:
\lim_{n\to\infty}
||S_{W_1\dots W_n}f -f^*||_\infty > 0 \right)
\hspace{-4cm}&&\\
&\le& P\left (\exists g\in\mathcal G:
\lim_{n\to\infty}
||S_{W_1\dots W_n}g -g^*||_\infty > 0 \right)\\
&\le& \sum_{g\in \mathcal G} 
P\left(\lim_{\ n\to\infty} ||S_{W_1\dots W_n}g -g^*||_\infty>0\right)\\
&=&0\,,
\end{eqnarray*}
proving Eq.~(\ref{conclusion-46a}).
\end{proof}

For the proof of Theorem~\ref{thm:iid}, we need one more lemma.

\begin{lem}[Identification of symmetric decreasing functions] 
\label{lem:OG} Let $f\in\Cc(\XX)$.

\begin{itemize}
\item[(a)] {\rm (by polarization).}
Let $W$ be a random variable on $\Omega$ 
whose distribution satisfies $\mu(R=0)=0$.
If the orbit of each $x\in\SS^{d-1}$
under $G=\{u\in\SS^{d-1}: (0,u)\in \supp \mu\}$
is dense in $\SS^{d-1}$, then
$$
\EE(\I(S_Wf))=\I(f) \quad \Longleftrightarrow\quad f=f^*\,.
$$

\item[(b)] {\rm (by Steiner symmetrization).}
Let $U$ be a random variable on $\SS^{d-1}$, 
and let $\mu$ be its probability distribution.
If the support of $\mu$
contains a generating set of directions for $O(d)$, then
$$
\EE(\I(S_Uf))=\I(f)
\quad \Longleftrightarrow \quad f=f^*\,.
$$
\end{itemize}
\end{lem}

\begin{proof} 
For part (a), suppose that 
$\EE(\I(S_Wf))=\I(f)$.
It follows from Lemma~\ref{lem:HL} that
$\I(S_\omega f)=\I(f)$, and hence $S_\omega f=f$,
for $\mu$-a.e. $\omega$.
This means that $f(\tau_\omega x)\ge f(x)$ for $\mu$-a.e. $\omega$
and all $x\in \XX$. 

Let $u\in G$.  By assumption, $\mu$ assigns strictly positive 
measure to each neighborhood of $(0,u)$ in $\Omega$.
Since $\mu(\{R=0\})=0$,
we can find a sequence $\omega_i=(r_i,u_i)$ with 
$r_i>0$ that converges to $(0,u)$
such that $f(\tau_{\omega_i}x)\ge f(x)$
for each $i$ and all $x\in \XX$.
By continuity, $f(\tau_u x)\ge f(x)$ for all $x\in\XX$,
which means that the value of $f$ increases
monotonically along orbits
$\tau_{u_n}\dots \tau_{u_1} x$ of $G$.
Since $\mathcal{O}_{G,x}$
is dense in the sphere
of radius $|x|$ and $f$ is uniformly continuous, $f$ must 
be radial.

To see that $f$ is symmetric decreasing, we write
it as $f(x)=\phi(d(x,o))$ for some continuous function $\phi$.
Consider first the cases $\XX=\RR^d$ and $\HH^d$.
Given $t>0$, choose $\omega=(r,u)$
with $0<r\le 2t$ such that 
$f(\tau_\omega x)\ge f(x)$
for all $x\in \XX$, and let $a$ be the point with
normal coordinates $(r,u)$.
The reflection $\sigma_{\omega}$ maps the centered
sphere of radius $t$ to the sphere of the same radius
centered at $a$. Since this sphere contains 
the points with normal coordinates $(r\pm t,u)$,
by the intermediate value theorem it 
contains for each $s\in (t,t+r]$ a
point $x$ with $d(x,o)=s$.
Since $d(\sigma_\omega x,o)=t<s$, the point
$x$ lies in the negative half-space $H^-_\omega$. 
It follows that
$\phi(s)=f(x)\le f(\tau_{\omega}(x)) = \phi(t) $.
Iterating the  argument, we conclude that
$\phi(s)\le \phi(t)$ for all $s\ge t>0$.
Since $\phi$ is continuous, we can take $t\to 0$ 
and conclude that $\phi$ is non-increasing on $[0,\infty)$.
In the case $\XX=\SS^d$, the above argument remains valid
for $t\in (0,\pi)$, provided that $r\le \min\{2t,\pi-t\}$,
and we obtain that $\phi$ is nonincreasing on $[0,\pi]$.
This proves that $f=f^*$. 

For part (b), suppose that $\EE(\I(S_{U}f))=\I(f)$.
We augment the random direction 
$U$ to a random variable $W=(R,\pm U)$
on $\Omega$, where $R$ is exponentially distributed on $\RR^+$,
the positive and negative signs are equally likely,
and the three components are independent.
Then $\EE(\I(S_{(R,\pm U)}f)=\I(f)$
by Eq.~\eqref{eq:trump}.
The probability distribution of $W$ is given by 
the measure $d\nu(r,u)= \frac{1}{2} e^{-r} dr (d\mu(u)+d\mu(-u))$
on $\Omega$.
By construction, $\nu(\{R=0\})=0$.
Since the support of $\mu$ contains a 
generating set of directions for $O(d)$,
the orbit
of any vector $x\in\SS^{d-1}$ under 
$$G= \bigl\{(0,u)\in \supp \nu\bigr\}= 
\bigl\{ \pm u: u\in\supp\mu\bigr\}$$
is dense in $\SS^{d-1}$.
Therefore, $\nu$ satisfies the assumptions of part (a),
and we conclude that $f=f^*$. Finally, the
converse implications hold 
because $S_\omega f^*=f^*$ for all $\omega\in\Omega$.
\end{proof}


\begin{proof}[Proof of Theorem~\ref{thm:iid}]
Let $W$ be a random variable on $\Omega$ that is distributed
according to the measure $\mu$ from the statement of the
theorem.  Lemma~\ref{lem:OG} guarantees that
$\EE(\I(S_Wf))<\I(f)$ unless $f=f^*$.
Let $\C_{L,\eta}$ be the set of all nonnegative continuous
functions supported in the ball of radius $L$ 
whose modulus of continuity is bounded by $\eta$.
Since $\I$ is continuous in the uniform topology
and $\C_{L,\eta}$ is compact by the the Arzel\`a-Ascoli theorem,
$$
h(\eps):= \inf \bigl\{ \EE(\I(f)-\I(S_Wf))
:  f\in \C_{L,\eta}, ||f-f^*||_\infty\ge\eps \bigr\} >0
$$
for each $\eps>0$. 

Given $f\in\Cc$, let $\eta$ be its modulus of continuity, and
assume that $f$ is supported in $B_L$.
Denote by $F_n=S_{W_1\dots W_n}f$
the result of $n$ random polarizations of $f$. Since 
polarization preserves the modulus of continuity and
the ball $B_L$, we have $F_n\in \C_{L,\eta}$.  
We argue as in the proof of Theorem~\ref{thm:46} that
\begin{eqnarray} 
\notag \I(f)-\I(f^*)
\notag &\ge& \sum_{n=1}^\infty \ 
\EE( \EE(\I(F_{n-1})-\I(S_{W_n}F_{n-1})
\ \vert\ W_1\dots W_{n-1})) \\
\label{eq:tatata}
 &\ge& h(\eps)\cdot \sum_{n=1}^\infty 
P( ||F_{n-1}\!-\!f^*||_\infty\ge \eps)\,.
\end{eqnarray}
In the second line, we have used the Markov property
and the definition of $h(\eps)$.
Since $h(\eps)>0$, the sequence
$F_n$ converges almost surely uniformly to $f^*$,
and Eq.~(\ref{conclusion-46a}) follows.  
\end{proof}


\begin{proof}[Proof of Corollary~\ref{cor:46}]
We proceed as
in the proofs of Theorems~\ref{thm:46} and~\ref{thm:iid},
with Eq.~\eqref{eq:key-Steiner} 
and Lemma~\ref{lem:OG}b in place of Eq.~\eqref{eq:key-pol}
and Lemma~\ref{lem:OG}a.
\end{proof}

\section{Examples in $\RR^d$}

The following lemma allows to transform integrals 
over $\Omega$ into integrals over $\RR^d$.
Geometrically, we map $\omega=(r,u)$  to the image of a point
$a$ under the reflection $\sigma_\omega$. Since for every
point $z\ne a$ there exists a unique reflection
that maps $a$ to $z$, this defines a diffeomorphism 
from $\Omega\setminus\{r=0\}$ to $\RR^d\setminus\{a\}$. 
For $a=o$, the diffeomorphism agrees with the polar coordinate map.

\begin{lem}[Change of variables] \label{lem:jacobian}
Let $a\in\RR^d$. 
Then 
$$
\int_\Omega g (\sigma_\omega a) \, d\omega
= \int_{\RR^d} g(z) |z-a|^{-(d-1)}\, dz\,
$$
for every measurable function $g$ on $\RR^d$ such
that the integral on the left hand side converges.
Here $d\omega=drdm(u)$ denotes the uniform measure on 
$\Omega$. 
\end{lem}

\begin{figure}[t]
\centerline{\includegraphics[width=0.36\textwidth]{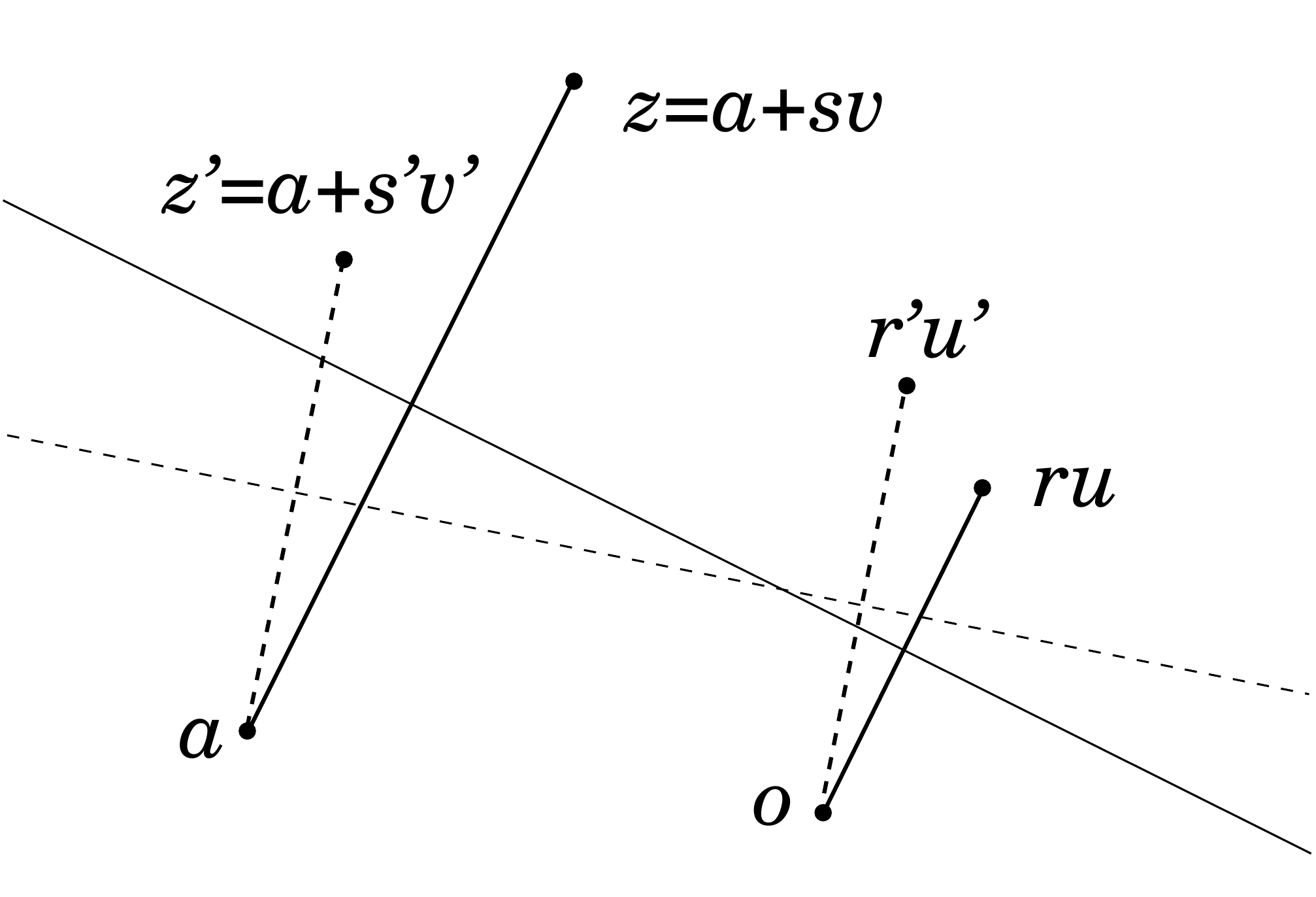}}
\caption{The change of variables in Lemma~\ref{lem:jacobian}.
In polar coordinates centered at 
$a$ and $o$,
the volume element transforms as $dsdv=drdu$.
}
\end{figure}

\begin{proof}
Set $z=\sigma_\omega a$.
If we write $\omega=(r,u)$ and express $z-a$ in polar coordinates $(s,v)$,
then $v=u$ because the
lines $x=\xi+tu$ are invariant under $\sigma_\omega$. 
If $r$ moves by a certain distance,
then $z$ moves by that distance in either the direction of $u$
or in the opposite direction (see Fig.~4).  In polar coordinates, the
metric on $\Omega$ transforms as $(ds)^2+(dv)^2=(dr)^2 +(du)^2$.
The claim follows by returning to Cartesian coordinates for $z$.
\end{proof}

We use this formula to construct examples of measures that 
satisfy the hypothesis of Theorem~\ref{thm:46} but
not Eq.~(\ref{assumption-VS}).  Consider the 
Gaussian probability measure on $\RR^d$
whose density is the centered heat kernel 
at time~$t$. By changing to polar coordinates,
we obtain a probability measure on $\Omega$, given by
$$
\mu(A)=\frac{1}{(2\pi t)^{\frac{d}{2}} }
\int_A e^{-\frac{r^2}{2t}} r^{d-1}\, d\omega\,,
$$
where $\omega=(r,u)$.
Fix $\rho, L>0$, let $a,b$ be a pair of points in $\RR^d$ with
$|a|+2\rho\le |b|\le L-\rho$, and consider the event
$\{\omega: d(\sigma_\omega a,b)<\rho\}$.
If $z=\sigma_\omega a\in B_\rho(b)$,
we use that $|z|-|a|<r< |z|+|a|$ to see
that $r\in [2\rho, 2L]$.  It follows that there
exists a constant $C$ (depending on $\rho$, $L$, and the dimension but
not on $t$) such that the density of $\mu$ in this region
is bounded from below by
$Ct^{-\frac{d}{2}} e^{-\frac{2L^2}{t}}$.
Changing variables with Lemma~\ref{lem:jacobian},
we estimate
$$
\mu(\{\omega:  d(\sigma_\omega a,b)<\rho\})
\ge  
Ct^{-\frac{d}{2}} e^{-\frac{2L^2}{t}}\int_{B_\rho(b)} |z-a|^{-(d-1)}\, dz\
\ge C't^{-\frac{d}{2}} e^{-\frac{2L^2}{t}}\,.
$$
Therefore $P$, the product of a sequence of such measures, satisfies
$$
\sum_{i=1}^\infty P(d(\sigma_{W_i} a_i,b_i)<\rho) \ge
C'\sum_{i=1}^\infty  t_i^{-\frac{d}{2}}
e^{-\frac{2L^2}{ t_i}}
$$
for any pair of sequences $\{a_i\}$, $\{b_i\}$ in $\RR^d$
with $|a_i|+2\rho \le |b_i|\le L-\rho$.
For $t_i=(\log\log i)^{-1}$ 
the sum diverges as required by Eq.~(\ref{assumption-46a}),
but Eq.~(\ref{assumption-VS}) fails because the measures converge
weakly to zero.  For the sequence $t_i = i^{2/d}$,
Eq.~(\ref{assumption-46a}) holds
but Eq.~(\ref{assumption-VS}) fails because
the measures concentrate on $\{0\}\times\SS^{d-1}$.

To give a similar example for Steiner symmetrizations, 
consider the probability measures on $\SS^{d-1}$
defined by the Poisson kernel
$$
\mu(A)=\frac{1}{m(\SS^{d-1})} \int_A
\frac{1-|z|^2}{|z-u|^{d}}\, dm(u) \,,
$$
where $z$ is a point in the ball, 
and $dm$ denotes integration with respect
to the standard Riemannian volume in $\SS^{d-1}$.
Since the  density of $\mu$ with respect to the uniform 
probability measure on $\SS^{d-1}$ is bounded from 
below by $2^{-(d-1)} (1-|z|)$, the product of such measures satisfies
$$
\sum_{i=1}^{\infty}P (d(U_i,v_i)<\rho)
\ge  2^{-(d-1)}\frac{m(B_\rho)}{m(\SS^{d-1})} 
\sum_{i=1}^\infty(1-|z_i|) 
$$
for every sequence $\{v_i\}$ in $\SS^{d-1}$.
If $z_i= (1-1/i)\,u$ for some $u\in \SS^{d-1}$,
then the sum diverges and Eq.~(\ref{assumption-46b}) holds,
but condition (\ref{assumption-VS}) fails because
the measures converge weakly to the point mass at $u$.

In principle, the proofs of Theorems~\ref{thm:46} and~\ref{thm:iid} 
imply weak-type bounds on the rate of convergence. 
Eq.~(\ref{eq:tata}) yields that
$$
P(||F_n\!-\!f^*||_\infty\ge \eps)
\le \frac{ \I(f)-\I(f^*)}{
C_\eps\cdot \sum_{i=1}^{n} \inf_{x,b} 
P(d(\sigma_{W_i}x,b)>\rho)}\,,
$$
where  $C_\eps$ and $\rho$ depend
on the modulus of continuity of $f$.
Similarly, since
$||F_n-f^*||_\infty$ is non-increasing,
Eq.~(\ref{eq:tatata}) yields that
$$
P(||F_n\!-\!f^*||_\infty\ge \eps)
\le \frac{ \I(f)-\I(f^*)}{
h(\eps)}\, n^{-1}\,,
$$
where $h(\eps)$  depends on the distribution of the
random polarizations and the modulus of continuity 
of $f$. For i.i.d.~uniform sequences of rearrangements,
we have a more explicit bound:

\begin{prop} [Rate of convergence for random 
polarizations] \label{prop:uniform}
If $\{W_i\}_{i\ge 1}$ is a sequence of 
independent uniformly distributed random variables 
on $(0,2L)\times\SS^{d-1}$, then
\begin{equation} \label{eq:rate}
\EE(||S_{W_1\dots W_n}f-f^*||_1)
\le 2d\, m(B_{2L})\, ||f||_\infty \, n^{-1}
\end{equation}
for every nonnegative bounded measurable 
function on $\RR^d$ with support in $B_L$. 
If, additionally, $f$ is H\"older continuous 
with modulus of continuity $\eta(\delta)\le c\delta^\alpha$
for some $\alpha\in (0,1]$ and $c>0$, then
\begin{equation} \label{eq:rate-L}
\EE(||S_{W_1\dots W_n}f-f^*||_\infty)
\le 10 cL^\alpha n^{-\frac{\alpha}{d+\alpha}}\,.
\end{equation}
\end{prop}

In the proof of the proposition, we will use the following lemma.

\begin{lem} [Expected drop in symmetric difference] 
\label{lem:uniform} 
If $W$ is a uniformly distributed random variable
on $(0,2L)\times\SS^{d-1}$, then 
$$
m(A\bigtriangleup A^*) -\EE(m(S_W A \bigtriangleup A^*))
\ge 
\frac{1}{2d\, m(B_{2L})}\, (m(A\bigtriangleup A^*))^2
$$
for every measurable set $A\subset B_L$ in $\RR^d$.
\end{lem}

\begin{proof} Fix $\omega\in\Omega$, and let
$H^+$ and $H^-$ be the half-spaces associated with $\omega$.
By construction, polarization swaps the portion of $A\setminus \sigma_\omega A$
that lies in $H^-$ with its mirror image in $H^+$
(see Fig.~1a).  Of these sets, precisely the portion 
of $A\setminus A^*$ whose reflection lies in $A^*\setminus A$
contributes to the symmetric difference
$A\bigtriangleup A^*$, twice, see Fig.~3b.  But this just means that 
\begin{equation}
\label{eq:Delta}
m(A \bigtriangleup A^*) -m(S_\omega A \bigtriangleup A^*) = 2m
(\{x\in A^*\setminus A: \sigma_\omega(x) \in A\setminus A^*\})\,.
\end{equation}
We compute the expectation, using Fubini's theorem and 
the change of variables from Lemma~\ref{lem:jacobian}.
The result is
\begin{eqnarray*}
m( A\bigtriangleup A^*)-
\EE(m(S_W A\bigtriangleup A^*)) \hspace{-4cm} &&\\
&=& 2\int_{A^*\setminus A} P(\sigma_W(x)\in A\setminus A^*)\, dx\\
&=& \frac{1}{L\, m(\SS^{d-1})} \int_{A^*\setminus A} \int_{A\setminus A^*} 
|x-z|^{-(d-1)}\, dzdx\\
&\ge & \frac{1}{C}
(m(A\bigtriangleup A^*))^2\,,
\end{eqnarray*}
where $C= 2d\,m(B_{2L})$. 
In the last step, we have used that the distance between
$x$ and $z$ is at most $2L$, and that $A$ and $A^*$
have the same volume.
Note that the Riemannian volume
of the unit sphere in $\RR^d$ is related to the
Lebesgue measure of the unit ball by $m(\SS^{d-1})=d\,m(B_1)$.
\end{proof}


\begin{proof}[Proof of Proposition~\ref{prop:uniform}]
Consider first the case where $f=\Chi_A$ for some
measurable set $A\subset B_L$, and let $A_n=S_{W_n\dots W_1}A$.
By Lemma~\ref{lem:uniform}, the Markov property,
and Jensen's inequality,
\begin{eqnarray*}
\EE(m(A_{n-1}\bigtriangleup A^*)) -\EE(m(A_n\bigtriangleup A^*))
\hspace{-5cm}&&\\
&=& \EE\bigl(
m(A_{n-1}\bigtriangleup A^*)-\EE(S_{W_n} A_{n-1}\bigtriangleup A^*)
\ \vert \ W_1,\dots,W_{n-1})\bigr)\\
&\ge & \frac{1}{C} \EE\bigl((m(A_{n-1}\bigtriangleup A^*))^2\bigr)\\
&\ge & \frac{1}{C} \bigl(\EE(m(A_{n-1}\bigtriangleup A^*))\bigr)^2\,,
\end{eqnarray*}
where $C=2d\,m(B_{2L})$.
This shows that $z_n=C^{-1}\EE(m(A_n\bigtriangleup A^*))$
satisfies the recursion relation $z_n\le z_{n-1}(1\!-\!z_{n-1})$.
Since $z_{n}^{-1}\ge z_{n-1}^{-1} + 1$ and 
$z_0^{-1}\ge d 2^{d+1}$, it follows that 
\begin{equation}
\label{eq:rate-A}
\EE(m(A_n\bigtriangleup A^*))
\ \le C(n+d 2^{d+1})^{-1} \,.
\end{equation}
If $f$ is a nonnegative bounded measurable function on $B_L$, 
we use the layer-cake principle to write
$$
||f-f^*||_1=\int_0^{\infty} m(\{f>s\}\bigtriangleup\{f^*>s\})\, ds\,,
$$
and likewise for $F_n= S_{W_n\dots W_n}f$. 
Since $f$ is bounded, the integrand
vanishes for $s>||f||_\infty$, and we obtain
from Eq.~\eqref{eq:rate-A}  that
$$
\EE(||F_n-f^*||_1) \le C ||f||_\infty(n+d 2^{d+1})^{-1}\,,
$$
proving the first claim.

If $f$ is H\"older continuous,
then $F_n$ and and $f^*$ are H\"older continuous with the
same modulus of continuity. 
Let $\eps=||F_n-f^*||_\infty$,
and set $\rho=\eta^{-1}(\eps/4)$.
Since $F_n $ differs from $f^*$ 
by at least $\eps/2$ on some ball
of radius $\rho$, we have
$||F_n -f^*||_{L_1} \ge \eps m(B_\rho)/2$.
We obtain from Eq.~(\ref{eq:rate}) that
\begin{eqnarray*}
\EE\bigl((||F_n-f^*||_\infty)^{1+\frac{d}{\alpha}}\bigr) 
&\le& \frac{2(4c)^{\frac{d}{\alpha}}}{m(B_1)}
\EE(||F_n-f^*||_1)\\
&\le& d 2^d (4cL^\alpha)^{1+\frac{d}{\alpha}} 
n^{-1}\,.
\end{eqnarray*}
Applying Jensen's inequality once more, we arrive at
$$
\EE(||F_n-f^*||_\infty) 
\le (d2^d)^{\frac{\alpha}{d+\alpha}}
4c L^\alpha n^{-\frac{\alpha}{d+\alpha}}\,.
$$
The leading constant is maximized 
at $\alpha=1$ and $d=6$, 
and Eq.~\eqref{eq:rate-L} follows.
\end{proof}

By Eq.~(\ref{eq:trump}), Proposition~\ref{prop:uniform}
extends to Steiner symmetrization along directions 
chosen independently and uniformly at random
on $\SS^{d-1}$.

\begin{cor} [Rate of convergence for random 
Steiner symmetrizations] \label{cor:uniform-S}
If $\{U_i\}_{i\ge 1}$ is a sequence of 
independent uniformly distributed random variables 
on $\SS^{d-1}$, then
$$
\EE(||S_{U_1\dots U_n}f-f^*||_1)
\le 2d\, m(B_{2L})\, ||f||_\infty n^{-1}
$$
for every nonnegative bounded measurable 
function with support in $B_L$.  
If $f$ is H\"older continuous with modulus of continuity 
$\eta(\delta)\le c\delta^\alpha$ for some $\alpha\in (0,1]$
and $c>0$, then
$$
\EE(||S_{W_1\dots W_n}f-f^*||_\infty)
\le 10 cL^\alpha n^{-\frac{\alpha}{d+\alpha}}\,.
$$
\end{cor}

\section{Negative results}

In this section, we give some bounds on the 
rate of convergence that complement Proposition~\ref{prop:uniform}
and Corollary~\ref{cor:uniform-S}, and construct examples 
where convergence fails.  For polarization, we use the function 
\begin{equation}
\label{eq:def-f-P}
f(x)=[1-d(x,a)]^+\,,
\end{equation}
which  is supported on $B_1(a)$ and
Lipschitz continuous with constant one.
Its symmetric decreasing rearrangement is
given by $f^*(x)=[1-|x|]^+$,
and $||f-f^*||_\infty =\min\{d(a,o),1\}$.
Its polarization at $\omega\in\Omega$ is given by
$$
S_\omega f(x)=[1-d(x,\tau_\omega a)]^+\,,
$$
where $\tau_\omega $ is the folding map that fixes the
positive half-space $H^+_\omega$ and reflects 
$H^-_\omega$ across the separating hyperplane.


\begin{prop} \label{prop:lower-P}
Let $f$ be given by Eq.~(\ref{eq:def-f-P}).

\begin{itemize}
\item[(a)]
{\rm (Convergence of random
polarizations is not faster than 
exponential).}\\
If $d(a,o)\le 1$, then
$$
\EE(||S_{W_1\dots W_n}f-f^*||_\infty\,) 
\ge ||f-f^*||_\infty\,2^{-n} \,
$$
for every sequence $\{W_i\}_{i\ge 1}$ of independent random variables
on $\Omega$ such that the distribution of
each $W_i=(R_i, U_i)$ is symmetric under $U_i\mapsto -U_i$. 

\item[(b)]  {\rm (Non-convergence).}
If $a\ne o$, then there exists a dense sequence 
$\{\omega_i\}_{i\ge 1}$ in $\Omega$ such that 
$S_{\omega_1\ldots \omega_n} f$ 
has no limit in $\Cc(\XX)$.
\end{itemize}
\end{prop}

\begin{proof}  A single random polarization
results in $S_Wf(x)= [1-d(x,\tau_Wa)]^+$.
Since $\tau_Wa=a$ whenever $a\in H_W^-$,
its expected distance from the origin satisfies
$\EE(d(\tau_Wa,o)) \ge d(a,o)/2$.  By iteration, we have
$S_{W_1\dots W_n} f(x)=[1-d(x,a_n)]^+$, 
where $a_n=\tau_{W_n}a_{n-1}$ and $a_0=a$.
By the Markov property, $ \EE(d(a_n,o)) \ge d(a,o)2^{-n}$,
and the first claim follows.

For the second claim, we realize an arbitrary sequence 
as a subsequence of one for which convergence fails.
Given $\{\omega_i\}_{i\ge 1}$,
fix $0<\eps< d(a,o)$ and define $\{\tilde \omega_i\}_{i\ge 1}$ as follows.
On the odd integers set
$\tilde \omega_{2n-1}=(\min\{2^{-n}\eps,r_n\},\pm u_n)$, where 
$(r_n,u_n)=\omega_n$,
and the sign is chosen in such a way that 
$S_{ \tilde \omega_1\dots \tilde \omega_{2n-1}}f$
is unchanged by $S_{\omega_n}$.
On the even integers, set $\tilde \omega_{2n}=\omega_n$.
If $\{\omega_i\}$ is dense, then $\{\tilde\omega_i\}$ is dense as well.

Set $f_n= S_{\tilde \omega_1\dots \tilde\omega_n}f=[1-d(x, a_n)]^+$. 
Suppose that $f_n$ converges to some limit~$g$. 
Then $g(x)=[1-d(x,b)]^+$ for some $b$.
Let $\omega=(r,u)\in\Omega$ with $r>0$. 
By density, we can find a subsequence 
$\{\tilde \omega_{n_k}\}$ that converges to $\omega$. 
Since both $a_{n_k-1}$ and $a_{n_k}=\tau_{\tilde \omega_{n_k}} a_{n_k-1}$
converge to $b$, we must have $\tau_\omega b=b$. 
Since $\omega$ was arbitrary, it follows that $b=o$.
On the other hand, $d(b,o)\ge d(a,o)-\eps\sum 2^{-n}>0$,
a contradiction.
\end{proof}

The corresponding bounds for Steiner symmetrizations on 
$\RR^d$ are slightly more involved.  As an example, we use the function 
\begin{equation}
\label{eq:def-f-S}
f(x)=[1-\langle x, Mx\rangle]^+\,,
\end{equation}
where $M$ is a positive definite symmetric $d\times d$ matrix.
The symmetric decreasing rearrangement of $f$
is $f^*(x)=[1-\lambda^* |x|^2]^+$, where
$\lambda^*$ is the geometric  mean of the eigenvalues 
of $M$.  The distance from $f$ to $f^*$ satisfies
\begin{equation}\label{eq:uniform-lambdas}
\frac{\lambda_{\max}-\lambda_{\min}}{2\lambda_{\max}}
\le ||f-f^*||_\infty\le 
\frac{\lambda_{\max}-\lambda_{\min}}{\lambda_{\min}}\,.
\end{equation}
We will prove the following statements.

\begin{prop} \label{prop:lower-S}
Let $f$ be given by Eq.~(\ref{eq:def-f-S})
with some positive definite symmetric matrix $M$.

\begin{itemize} 
\item[(a)] {\rm  
(Convergence of random Steiner symmetrizations is not faster than 
exponential).}
If $\{U_i\}$ is a sequence of i.i.d.~uniform random
variables on $\SS^{d-1}$ and 
the extremal eigenvalues
of $M$ satisfy $\lambda_{\max} \le 2\lambda_{\min}$, then
$$
\EE(||S_{U_1\dots U_n}f-f^*||_\infty)
\ge \frac{1}{4} ||f-f^*||_\infty \, 3^{-n}\,.
$$

\item[(b)] {\rm (Non-convergence).}
If $M$ is not a multiple of the 
identity, then there exists a 
dense sequence $\{u_i\}_{i\ge 1}$ in $\SS^{d-1}$ such that
$S_{u_1\ldots u_n} f$ has no limit in $\Cc(\RR^d)$.
\end{itemize}
\end{prop}

We first show that Steiner symmetrization preserves the form of $f$.

\begin{lem} [Steiner symmetrization of ellipsoids]
\label{lem:ellipsoid}
If $f$ is given by
Eq.~(\ref{eq:def-f-S}), then $S_uf$ has the same
form with a positive definite symmetric matrix $M'$
determined by
\begin{equation} 
\label{eq:ellipsoid}
\ \langle x, M' x\rangle  = \langle x, M x\rangle 
-\frac{\,\langle x, Mu\rangle^2}{\langle u, Mu\rangle}
+ \langle x, u\rangle^2 \langle u,  Mu\rangle\,.
\end{equation}
In particular, $u$ is an eigenvector of $M'$ with
eigenvalue $\langle u, Mu\rangle$.
\end{lem}

\begin{proof} Consider a line $x=\xi+tu$ with $\xi\perp u$.
The restriction of $f$ to this line, given by
$$
t\mapsto 
[1-\langle \xi, M\xi\rangle  - 2t \langle \xi, Mu\rangle - 
t^2 \langle u, Mu\rangle ]^+\,,
$$
is symmetric decreasing about 
$t_0=-\frac{\langle \xi, Mu\rangle }{\langle u, Mu\rangle }$.  
By definition, the restriction of $S_uf$ to the
line is the symmetrized function
$$
t\mapsto 
[ 1- \langle \xi, M\xi\rangle  + (t_0^2-t^2)\langle u, M u\rangle ]^+\,,
$$
as required by Eq.~(\ref{eq:ellipsoid}). 
Since $u^\perp$ and the line through $u$
are invariant subspaces for $M'$, we conclude that $u$ is an
eigenvector.
The corresponding eigenvalue is $\lambda=\langle u,M'u\rangle=
\langle u,Mu\rangle$.
\end{proof}


\noindent{\em Remark.} 
An amusing consequence of Lemma~\ref{lem:ellipsoid}
is that $(d\!-\!1)$ Steiner symmetrizations suffice to
transform an ellipsoid into a ball~\cite{KM}.  
To see this let $A$ be an ellipsoid of the same volume as the unit
ball.
Then $A=\{\langle x, Mx\rangle<1\}=\{f>0\}$, 
where $M$ is a  positive definite symmetric
matrix of determinant one, and $f$ is given by
Eq.~(\ref{eq:def-f-S}).
Set $M_0=M$, and choose $u_1$ 
such that $\langle u_1, Mu_1\rangle =1$. By Lemma~\ref{lem:ellipsoid},
$S_{u_1}A= \{\langle x, M_1x\rangle <1\}$, where $M_1$ is a 
positive definite symmetric
matrix that has $u_1$ as an eigenvector
with eigenvalue~$1$.
Iteratively choosing $u_i$ orthogonal to $u_1,\dots, u_{i-1}$
such that $\langle u_i, M_{i-1} u_i\rangle =1$, we arrive at 
$M_{d-1}=I$, and conclude that $S_{u_{d-1}}\dots S_{u_1}A=A^*$.
\hfill $\Box$

\medskip 
To prove Proposition~\ref{prop:lower-S}, we need to 
analyze how the extremal eigenvalues of $M$ change under Steiner 
symmetrization of $f$.  Clearly, their difference decreases,
because the inradius of the corresponding
ellipsoid grows under Steiner symmetrization, 
and its outradius shrinks. The following lemma shows that the
change in the extremal eigenvalues is small, if the direction of
the Steiner symmetrization is either almost parallel
or almost orthogonal to the maximizing eigenvector 
$v_{\max}$ (see Fig.~5).

\begin{figure}[t]
\centerline{\includegraphics[width=0.32\textwidth]{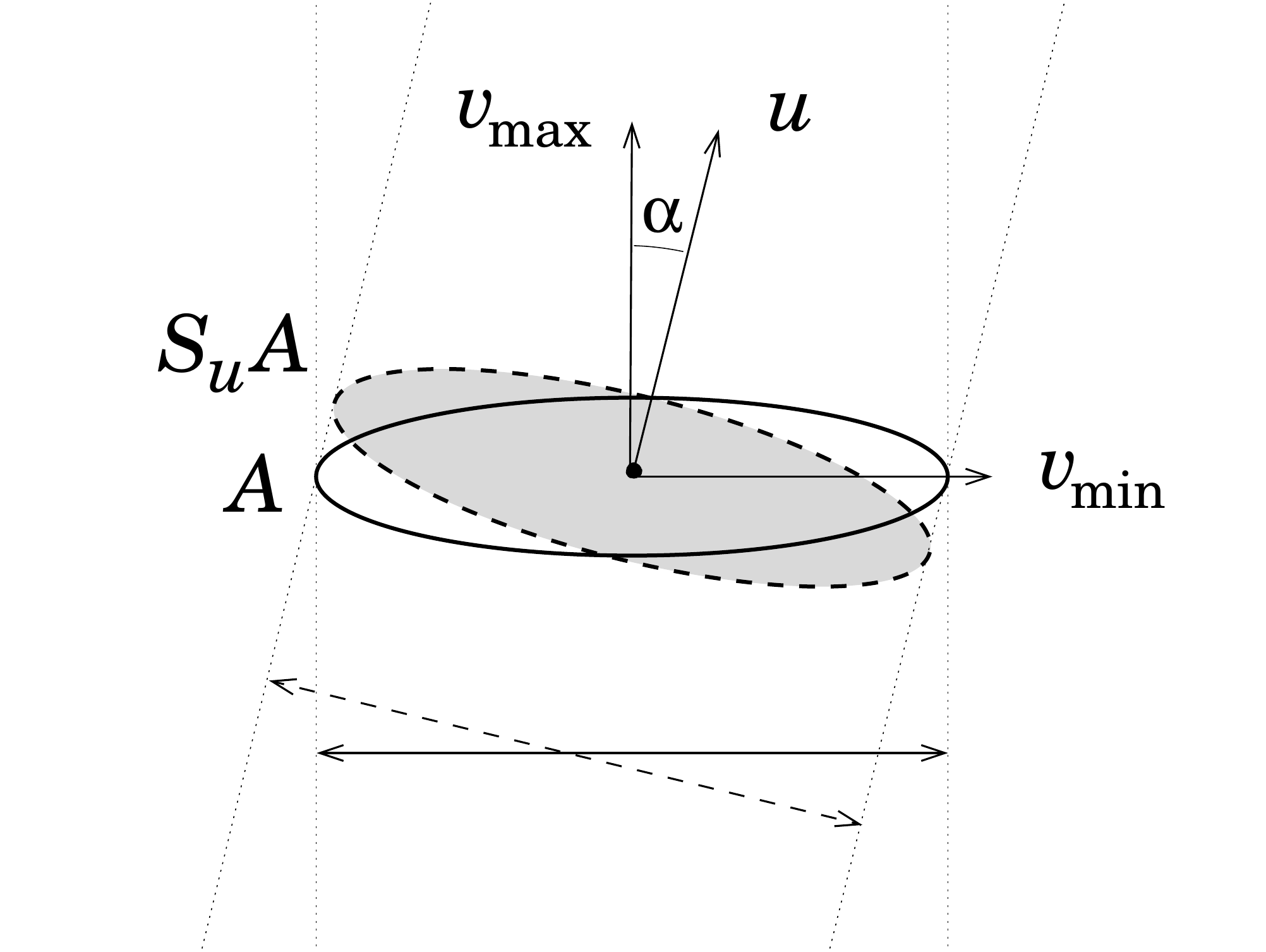}}
\caption{Steiner symmetrization of an ellipse.
The diameter shrinks at most by a factor $\cos\alpha$.
} 
\end{figure}

\begin{lem} [Eigenvalue estimate]
\label{lem:eval}
Given a symmetric positive definite matrix $M$
with extremal eigenvalues $\lambda_{\max}$, $\lambda_{\min}$ 
and corresponding normalized eigenvectors $v_{\max}$, $v_{\min}$.
Define $M'$ by Eq.~(\ref{eq:ellipsoid}).
The extremal eigenvalues $\lambda'_{\max}$, $\lambda'_{\min}$
of $M'$ satisfy
\begin{equation}
\label{eq:eval}
\lambda'_{\max}-\lambda'_{\min}
\ge 
\bigl(1-C\psi(\langle u,v_{\max}\rangle) -2 
\psi(\langle u, v_{\min}\rangle) \bigr)\,(\lambda_{\max}-\lambda_{\min})\,,
\end{equation}
where $C=1+\lambda_{\max}/\lambda_{\min}$ and
$\psi(t)=t^2(1-t^2)$.
\end{lem}

\begin{proof} Let $v$ be a normalized eigenvector of $M$ 
with eigenvalue $\lambda$. From Eq.~(\ref{eq:ellipsoid}), we obtain that
\begin{eqnarray*}
\langle v,  M' v\rangle
&=& \lambda - \frac{\lambda^2 \langle u, v\rangle ^2}{\langle u, Mu\rangle}
+ \langle u,v\rangle^2 \langle u,Mu\rangle\\
&=& \lambda + \cos^2\!\alpha \sin^2\!\alpha \,
\left(1+\frac{\lambda}{\langle u, Mu\rangle}\right)
(\langle w, Mw\rangle -\lambda)\,.
\end{eqnarray*}
In the second step, we have expanded 
$u=\cos\alpha \,v + \sin\alpha\,w$,
where $w$ is a unit vector orthogonal to $v$,
and then collected terms.
We apply this identity to $v_{\max}$
and use that $\langle u, Mu\rangle \ge \lambda_{\min}$ and 
$\langle w, Mw\rangle \le \lambda_{\max}$
to obtain
\begin{eqnarray*}
\lambda'_{\max}
&\ge &\langle v_{\max},M'v_{\max}\rangle\\
&\ge&
\lambda_{\max}- 
\left(1+\frac{\lambda_{\max}}{\lambda_{\min}}\right)\,
\psi(\langle u,v_{\max}\rangle) \, (\lambda_{\max}-\lambda_{\min})\,.
\end{eqnarray*}
Similarly, 
\begin{eqnarray*}
\lambda'_{\min}&\le& \langle v_{\min},M'v_{\min}\rangle\\
& \le& \lambda_{\min} + 
2 \psi(\langle u, v_{\min}\rangle) \,(\lambda_{\max}-\lambda_{\min})\,.
~~~~~~~~~~~~~~~
\end{eqnarray*}
The claim follows by subtracting the two inequalities.
\end{proof}


\begin{lem} [Expected change of extremal eigenvalues]
\label{lem:extremal}
Let $f$ be given by Eq.~(\ref{eq:def-f-S}) with
a positive definite symmetric matrix $M$
whose extremal eigenvalues 
satisfy $\lambda_{\max}\le 2\lambda_{\min}$.
If $U$ is a uniformly distributed random variable
on $\SS^{d-1}$, then $S_Uf(x)= [1-\langle x, M'x\rangle ]^+$, where
$M'$ is a positive definite symmetric matrix whose
the extremal eigenvalues satisfy
$$ \EE(\lambda'_{\max}-\lambda'_{\min}) \ge 
\frac13(\lambda_{\max}-\lambda_{\min}) \,.
$$
\end{lem}

\begin{proof} We apply Lemma~\ref{lem:eval} and take expectations.
Let $v_{max}$ and $v_{\min}$
be the eigenvectors of $M$ corresponding to $\lambda_{\max}$
and $\lambda_{\min}$, and set
$C=1+\lambda_{\max}/\lambda_{\min}\le 3$
and $\psi(t)=t^2(1-t^2)$.
By taking advantage of the rotation invariance,
we compute
$\EE(\langle U, v\rangle^2) = 1/d$
and $\EE(\langle U,v\rangle^4\,)
= 3/(d(d+2))$ for all $v\in \SS^{d-1}$,
see~\cite[Exercise 63, p. 80]{Folland}. 
This results in
$$\EE(1-C\psi(\langle U, v_{\max}\rangle -2\psi(\langle U, 
v_{\min}\rangle)) 
= 1-(C+2)\left(\frac{1}{d}-\frac{3}{d(d+2)}\right) \,.
$$ 
The claim follows by evaluating the right hand side
at $d=3$, where it assumes its minimum value,
and using Eq.~(\ref{eq:eval}).
\end{proof}


\begin{proof} [Proof of Proposition~\ref{prop:lower-S}]
Let $f$ be given by 
Eq.~(\ref{eq:def-f-S}) with some positive definite symmetric matrix $M$.
We first consider the case of a random sequence
$F_n=S_{U_1\dots U_n}f$, where the directions
$\{U_i\}$ are independent and uniformly distributed on $\SS^{d-1}$.
By Lemma~\ref{lem:ellipsoid}, we can write
$F_n$ in the form (\ref{eq:def-f-S}) with a positive definite
symmetric matrix $M_n$ 
that is recursively defined by Eq.~(\ref{eq:ellipsoid})
with $u=U_n$.  We iterate the estimate in Lemma~\ref{lem:extremal},
using the Markov property,
and obtain that the gap between the extremal eigenvalues of $M_n$
is at least $(\lambda_{\max}\!-\!\lambda_{\min}) 3^{-n}$.  
Since we assumed that $\lambda_{\max}\le 2\lambda_{\min}$,
it follows from Eq.~(\ref{eq:uniform-lambdas}) 
that 
\begin{eqnarray*}
\EE( ||S_{U_1\dots U_n}f-f^*||_\infty\,)
&\ge& \frac{\lambda_{\max}-\lambda_{\min}}{2\lambda_{\max}}3^{-n} \\
&\ge& \frac 14||f-f^*||_\infty \, 3^{-n}\,.
\end{eqnarray*}

For the second claim, we proceed as
in the proof of Proposition~\ref{prop:lower-P}
by realizing an arbitrary sequence 
as a subsequence of one for which convergence fails.
Given $\{u_i\}_{i\ge 1}$ in $\SS^{d-1}$, 
let $\eps>0$ so small that $ (C+2)\sin^2 \eps<1$,
where $C=1+\lambda_{\max}/\lambda_{\min}$ 
as in Lemma~\ref{lem:extremal},
and construct the sequence $\{v_i\}_{i\ge 1}$ as follows.  
In the first step, pick $v_1$
to be a maximizing eigenvector of $M$. 
Suppose we have already chosen 
$v_1,\dots ,v_n$ such that $d(v_i,v_{i+1})\le \eps/i$
for each $i<n$, and that $u_1,\dots ,u_j$
appear as a subsequence. If $d(v_n,u_{j+1})\le \eps/n$, 
pick $v_{n+1}=u_{j+1}$.  Otherwise, choose $v_{n+1}$ on the 
great circle that joints $v_n$ with $u_{j+1}$
in such a way that
$d(v_n,v_{n+1})=\eps/n$ and 
$d(v_{n+1},u_{j+1}) =d(v_n,u_{j+1})- \eps/n$.
Since $\sum \eps/n$ diverges, 
the entire sequence $\{u_i\}$ is incorporated as
a subsequence into $\{v_i\}$. If $\{u_i\}$ is dense, so is $\{v_i\}$.  

Let $f_n=S_{v_1\dots v_n}f=[1-\langle x, M_n x\rangle]^+$.
If $f_n$ converges to some limit~$g$, then $g$ is given by 
Eq.~(\ref{eq:def-f-S}) with some
positive definite symmetric matrix $N$ 
of the same determinant as $M$.  Since 
$\{v_i\}$ is dense in $\SS^{d-1}$, 
we find that $N$ is necessarily 
a multiple of the identity because
$g$ is invariant under every Steiner symmetrization.
On the other hand, we can estimate
the extremal eigenvalues of $N$ as follows.
By construction, $v_n$ is an eigenvector
of $M_n$.  Since $d(v_n, v_{n+1})\le \eps/n$
and the other eigenvectors of $M_n$ are orthogonal to $v_n$,
we have that 
$\psi(\langle v_{n+1},v\rangle)\le \sin^2 (\eps/n)$
for each eigenvector $v$ of $M_n$. Iterating
Lemma~\ref{lem:eval}, we see that
the gap between the extremal eigenvalues of $N$
is at least
$(\lambda_{\max}-\lambda_{\min})\prod(1-(C+2)\sin^2(\eps/n))>0$, 
a contradiction.
\end{proof}

\section{Compact sets}

We finally collect the implication of our results for 
compact sets. 
Under the assumptions of Theorems~\ref{thm:46} and~\ref{thm:iid}, 
random polarizations of functions in $L_p^+$
also converge almost surely in $L_p$,
see Eq.~(\ref{eq:Lp}).  In particular, for $p=1$,
\begin{equation}\label{eq:symmdiff}
P\left(\,\lim_{n\rightarrow \infty}
m(S_{W_1\dots W_n}A\bigtriangleup A^*) =0 
   \quad \forall A\subset \XX\ \mbox{with}\  m(A)<\infty \right)=1 \,.
\end{equation}
This is another equivalent restatement of Eq.~(\ref{conclusion-46a}).
We now establish the corresponding convergence result
for the Hausdorff distance.

The topology defined by the Hausdorff metric on the space of compact sets
is not comparable to the topology of symmetric difference.
Moreover, polarization is not continuous with respect to Hausdorff distance.
To give a simple example,
consider a reflection $\sigma$ that does not fix the origin,
and let $a$ be the image of the origin under $\sigma$.
By definition, $S(\{o,a\})=\{o,a\}$. 
Let $\{a_i\}_{i\ge 1}$ be a sequence in $\XX$ 
with $a_i\ne a$ for all $i$ that converges to $a$.
The sequence of two-point sets $\{o,a_i\}$ clearly converges to
$\{o,a\}$. Since the reflected sequence $\{\sigma a_i\}$
converges to the origin, we have 
that $d(a_i,o)>d(\sigma a_i,o)$ and therefore $a_i\in H^-$
for $i$ large enough. It follows that
$S(\{o,a_i\})= \{o,\sigma a_i\}$,
which converges in Hausdorff distance to $\{o\}$.

Nevertheless, convergence of a sequence of polarizations  in Hausdorff 
distance  to a ball implies convergence in symmetric difference.
To see this, let $K$ be a compact set of positive volume,
and consider a sequence $K_n=S_{\omega_1\dots\omega_n}K$. 
If $K_n $ converges to (the closure of) $K^*$ in 
Hausdorff distance, then the radius of the smallest centered
ball containing $K$ converges to the radius of $K^*$,
which implies that $m(K_n \setminus K^*)$
converges to zero.  Since $K_n$ and $K^*$
have the same volume, $m(K^*\setminus K_n)$
goes to zero as well. 

In the other direction, we can obtain convergence 
in Hausdorff distance from the uniform convergence
statement in Eq.~\eqref{conclusion-46a}
by realizing a given compact set as a level set of a continuous function.

\begin{prop} [Convergence in Hausdorff distance]
\label{prop:compact}
If a random sequence $\{W_i\}$ 
satisfies the assumptions 
of either Theorem~\ref{thm:46} or Theorem~\ref{thm:iid}, then
$$
P\left( \lim_{n\to\infty} d_H\bigl(S_{W_1\dots W_n}K,K^*\bigr) =0
\quad \forall \ \mbox{compact}\ K\subset \XX 
\ \mbox{with}\ m(K)>0
\right) =1\,,
$$
and 
$$
P\left( \! \lim_{n\to\infty} d_H\bigl(\partial 
S_{W_1\dots W_n}K,\partial K^*\bigr) =0
\;
\begin{array}{l} \forall \ 
\mbox{compact}\ K\subset \XX \\
\mbox{with}\ m(K)>0\
\mbox{and}\ m(\partial K)=0 \end{array}\!\!\right) =1\,.
$$
\end{prop}

\smallskip

\begin{proof} Set $K_n=S_{W_1\dots W_n}K$.  We consider 
the two pieces of the Hausdorff distance from $K_n$ to $K^*$ separately.
If $\dist(x,K_n)=\delta>0$ for some $x\in K^*$, then
$m(K_n \bigtriangleup K^*) \ge 2m(B_\delta (x)\cap K^*)>0$.
Therefore Eq.~(\ref{eq:symmdiff}) implies that
$$
\sup_{x\in K^*} \dist(x,K_n) \to 0
\quad (n\to\infty)\quad\mbox{almost surely}
$$
simultaneously for all $K$.  

To control the other piece of $d_H(K_n,K^*)$, we use the auxiliary function
$ f(x)=[1 - \dist(x,K)]^+$.
By definition, the level set of $f$ at height $1-t$
is the outer parallel set $\{x:\dist(x,K)< t\}$.
The level set of $f^*$ at that height is the
centered ball of the same volume.
Its radius $\rho(t)$, defined by
$$
B_{\rho(t)} = \{x: \dist(x,K)<t\}^* \quad (t> 0)
$$
depends continuously on $t$
and converges to the radius of $K^*$ as $t\to 0$.
Set $F_n=S_{W_1\dots W_n}f$.  Since $F_n(x)=1$ for all $x\in K_n$,
\begin{eqnarray}
\notag
\sup_{x\in K_n} \dist(x, K^*) 
&=& \sup_{x\in K_n\setminus K^*} \rho(F_n(x)\!-\!f^*(x)) -
\mbox{radius}\,(K^*) \\
\label{eq:dist-out}
& \to& 0\quad (n\to\infty)\quad  \mbox{almost surely}
\end{eqnarray}
by Theorem~\ref{thm:46}.  This proves the first claim.

If $\partial K$ has zero volume, we continuously extend
the function $\rho$ such that 
$$
B_{\rho(0)}=K^*\,,\qquad B_{\rho(t)} = \{x: \dist(x,\XX\setminus K)> -t\}^* 
\quad (t<0)\,,
$$
and replace the auxiliary function with
\begin{equation} \label{eq:f-dist}
f(x)=\bigl[h+\dist(x,\XX\setminus K)-\dist(x,K)\bigr]^+\,, 
\end{equation}
where $h>0$ is an arbitrary constant.
The level sets of $f$ at heights below $h$ are outer parallel sets of $K$,
while the level sets at heights above $h$
are inner parallel sets.  It follows that
\begin{eqnarray}
\notag
d_H(\partial K_n,\partial K^*)
&=& 
\sup_{x\in \partial K_n} |\rho(h-f^*(x))-\mbox{radius}\,(K^*)|\\
\label{eq:radii-s}
&\le& \max_{\pm}
|\rho(\pm||F_n\!-\!f^*||_\infty)- 
\rho(0)| \\
\notag
&\to& 0 \quad (n\to\infty)\quad \mbox{almost surely.}
\end{eqnarray}
In the second line, we have used that
$F_n=h$ on $\partial K_n$.  The last line follows from
Theorem~\ref{thm:46} and the continuity of $\rho$.
\end{proof}

Similar arguments can be used to bound the rate
of convergence for sets with additional regularity 
properties.  Let $K$ be a compact set in $\RR^d$, and define
$f$ and $\rho$ as in the proof of
the second claim of Proposition~\ref{prop:compact}.
Assume that $K\subset B_L$, and that
$\rho$ is differentiable at $t=0$
with $\rho'(0) = \Per(K)/\Per(K^*)$.
By Proposition~\ref{prop:uniform} there exists
a sequence $\{\omega_i\}$ such that
$$
|| S_{\omega_1\dots\omega_n} f-f^*||_\infty
\le 10 (L+h)\, n^{-\frac{1}{d+1}}\,.
$$
Expanding $\rho$ about $t=0$, we obtain from Eq.~(\ref{eq:radii-s})  that
\begin{align}
\notag
d_H(\partial K_n,\partial K^*)
&\le 
\rho'(0)\, (1+o(1)) \, ||S_{\omega_n\dots \omega_1}f
-f^*||_\infty \\
\label{eq:radii}
&\le 
C\cdot{\rm radius}\,(K^*) \cdot \frac{\Per(K)}{~\Per(K^*)}\,
(1+o(1))\,n^{-\frac{1}{d+1}}
\end{align}
as $n\to\infty$, where $C=10 (L+h)/\rho(0)$.
After dropping an initial segment $n\le N$ from the sequence,
we may replace $L$ with the radius of the 
smallest centered ball containing $K_N$.
Choosing $N$ sufficiently large and $h$ sufficiently small, we can find 
a sequence of polarizations where Eq.~\eqref{eq:radii} holds with $C=10$.

\medskip
\noindent{\em Remark.}
The conclusions of Proposition~\ref{prop:compact}
also hold for random Steiner symmetrizations 
that satisfy the assumptions of Corollary~\ref{cor:46}.
Likewise, Eq.~\eqref{eq:radii} applies to
sequences of Steiner symmetrizations along i.i.d.~uniformly 
distributed directions.  However, in view of 
Klartag's result for convex sets, we expect 
such sequences to converge more rapidly (see Eq.~(\ref{eq:Klartag})).

\subsection*{Acknowledgments} 

This paper is based on 
results from M.F.'s 2010 Master's thesis
at the University of Toronto~\cite{Fortier}.
Our research was supported in part by NSERC 
through Discovery Grant No. 311685-10 (Burchard) and an 
Alexander Graham Bell Canada Graduate Scholarship (Fortier). 
A.B. wishes to thank Gerhard Huisken 
(Albert-Einstein Institut in Golm),
Bernold Fiedler (Freie Universit\"at Berlin),
Nicola Fusco (Universit\`a di Napoli Federico II) and Ad\`ele Ferone
(Seconda Universit\`a di Napoli, Caserta)
for hospitality during a sabbatical in 2008/09.
Special thanks go to Aljo\v{s}a Vol\v{c}i\v{c} for an
inspiring discussion of his results on Steiner symmetrization
that provided the original motivation for our work,
and to Bob Jerrard for pointing out an error in
an earlier version of Eq.~(\ref{eq:radii}).
We note that a non-convergence result
similar to Proposition~\ref{prop:lower-S}b was obtained 
independently by Bianchi, Klain, Lutwak, Yang, and
Zhang~\cite{BKLYZ}.


\baselinestretch\linespread{0.7} 

\providecommand{\bysame}{\leavevmode\hbox to3em{\hrulefill}\thinspace}
\providecommand{\MR}{\relax\ifhmode\unskip\space\fi MR }
\providecommand{\MRhref}[2]{%
  \href{http://www.ams.org/mathscinet-getitem?mr=#1}{#2}
}
\providecommand{\href}[2]{#2}

\end{document}